\def\AACA{0}  
\def\ELS{0}   
\def\RSOS{0}  
\def\ARXIV{0} 

\ifnum1=\AACA  
\documentclass[runningheads]{llncs}
\usepackage[T1]{fontenc}
\usepackage{graphicx} 
\usepackage{hyperref}
\usepackage{color}

\def\toobig{\vskip -0ex}
\else
\newcommand{\titlerunning}{}
\ifnum1=\RSOS
\documentclass{rsos} 
\titlehead{Algebra}
\def\toobig{\vskip -2ex}
\else
\ifnum1=\ELS
\documentclass{cas-sc}
\usepackage[numbers]{natbib}
\def\toobig{\vskip -0ex}
\newtheorem{theorem}{Theorem}
\newtheorem{lemma}{Lemma}
\newtheorem*{corollary}{Corollary}
\theoremstyle{definition}
\newtheorem*{definition}{Definition}
\else
\def\ARXIV{1} 
\documentclass[11pt]{article}
\usepackage[english]{babel}
\usepackage[a4paper,top=2cm,bottom=2.5cm]{geometry}
\usepackage[colorlinks=true, allcolors=blue]{hyperref}
\usepackage{graphicx}
\usepackage{amsthm}
\newtheorem{theorem}{Theorem}
\newtheorem{lemma}{Lemma}
\newtheorem*{corollary}{Corollary}
\theoremstyle{definition}
\newtheorem*{definition}{Definition}
\def\toobig{\vskip -2ex}
\fi\fi\fi

\usepackage{amsmath}
\usepackage{amssymb}
\usepackage{multirow}
\usepackage{anyfontsize}
\newcommand\Z{{\mathbb Z}} 
\newcommand\N{{\mathbb N}} 
\newcommand\R{{\mathbb R}} 
\newcommand\C{{\mathbb C}} 
\newcommand\HB{{\mathbb H}} 
\newcommand\OB{{\mathbb O}} 
\newcommand\SB{{\mathbb S}} 
\newcommand\PB[1]{{\mathbb P_{#1}}}

\newcommand\SCB{{\rm \widetilde{\mathbb C}}}
\newcommand\SHB{{\rm \widetilde{\mathbb H}}}
\newcommand\SOB{{\rm \widetilde{\mathbb O}}}
\newcommand\UB[1]{{\mathbb U_{#1}}}
\newcommand\vv[1]{\text{\bf #1}}
\newcommand\va[1]{\text{\bf a}_{#1}}
\newcommand\vb[1]{\text{\bf b}_{#1}} 

\newcommand\MA[1]{{\mathrm{A}(#1)}}
\newcommand\MB[1]{{\mathrm{B}(#1)}} 
\renewcommand\AA[1]{{\mathbb A}_{#1}}
\renewcommand\o[1]{{\rm o}_{#1}} 
\renewcommand\u[1]{{\rm u}_{#1}} 
\newcommand\e[1]{{\rm e}_{#1}} 
\newcommand\GT{{\rm G_2}}
\newcommand\coindex[1]{{\overline{#1}}}
\newcommand\salpha[2]{{\alpha_{#2}{#1}\hspace{1pt}
     \coindex{\alpha_{#2}}\hspace{1pt}}}
\newcommand\sbeta[2]{{\beta_{#2}{#1}\coindex{\beta_{#2}}\hspace{1pt}}}
\newcommand\sdelta[2]{{\delta_{#2}{#1}\coindex{\delta_{#2}}\hspace{1pt}}}

\newcommand\ec{\text{,}}
\newcommand\nz{$\star$}
\newcommand\es{\text{.}}

\newcommand\tstrut{\rule{0pt}{2.8ex}}
\newcommand\tab[3]{\begin{table}[#1]\toobig\centering\caption{\label{tab:#2}#3}}
\newcommand\tae{\end{table}}
\newcommand\tar[1]{Table~\ref{tab:#1}}
\newcommand\tas[1]{Tables~\ref{tab:#1}}
\DeclareMathSymbol{\sm}{\mathbin}{AMSa}{"39}
\def\abs{\mathop{\rm abs}\nolimits}

\newcommand\keywordlist{Octonions, sedenions, split seden\-ions, Cayley-Dickson algebras, zero divisors, graded algebras, Moufang loops, Malcev algebra, quasi-octonions, ultracomplex numbers}

\begin{document}
\ifnum1=\ELS
\title[mode=title]{Structure of the Cayley-Dickson algebras}
\shorttitle{Structure of the Cayley-Dickson algebras}
\tnotetext[1]{The document was supported by an Australian Government Research Training Program Scholarship and funding from the Australian Research Council (FL240100217).}
\tnotemark[]
\author{G.~P.~Wilmot}[orcid=0009-0003-7451-3933]
\shortauthors{G.~P.~Wilmot}
\ead{Greg.Wilmot@adelaide.edu.au}
\ead[url]{www.adelaide.edu.au}
\affiliation{organisation={School of EME, Adelaide University},
   city={Adelaide}, postcode={5005}, state={South Australia}, country={Australia}}
\else
\title{Structure of the Cayley-Dickson algebras}
\ifnum1=\RSOS
\author{G.~P.~Wilmot}
\address{University of Adelaide, Adelaide, South Australia, 5005, Australia, \url{https://www.adelaide.edu.au}}
\subject{Mathematics}
\corres{email:greg.wilmot@adelaide.edu.au}
\else
\ifnum1=\AACA
\titlerunning{Structure of the Cayley-Dickson algebras}
\author{G.\:P.\:Wilmot\inst{}\orcidID{0009-0003-7451-3933}}
\authorrunning{G.\:P.\:Wilmot}
\institute{Adelaide University, Adelaide, South Australia, 5005, Australia
\email{greg.wilmot@adelaide.edu.au}\\
\url{https://www.adelaide.edu.au}}
\else
\author{G.~P.~Wilmot \\
  University of Adelaide, Adelaide, South Australia, 5005, Australia,\\
  \url{https://www.adelaide.edu.au}
  email:greg.wilmot@adelaide.edu.au}
\fi
\thispagestyle{empty}
\maketitle
\fi\fi

\begin{abstract}
Viewing the Cayley-Dickson process as a graded construction provides a rigorous definition of associativity consisting of three classes, which divides the non-associ\-ative parts into four types. These simplify the Moufang loop identities and Malcev's identity and categorise non-associative loops. Analysing the non-associativity structure of sedenions and higher level Cayley-Dickson algebras uncovers octonion and three quasi-octonion, power-associative subalgebras. Power-associativity is proved to grow by multiples of seven that introduce zero divisors in factors of 84 for all Cayley-Dickson algebras. It is convenient to replace the terminology hypercomplex numbers with {\it ultracomplex numbers} to capture the fact that sedenions and all higher level Cayley-Dickson algebras contain quasi-octonions that are ultracomplex. The non-associative types reduce the 84 zero divisors in sedenions via three cycles and four modes down to seven primary zero divisor pairs of two-term sums of basis elements providing new techniques to analyse the seven quasi-octonions loops in sedenions. The graded notation allows the eight octonion and seven quasi-octonion subalgebras of sedenions to be uniquely derived, up to representation.

The zero divisors for split sedenion algebras are analysed and mappings between three of these are described showing 180 zero divisors and the other split algebra is analogous to sedenions with 84 zero divisors. These split algebras are shown to involve quasi-octonion algebras with mixed signature.
\end{abstract}
\ifnum1=\ELS
\maketitle
\else
\ifnum1=\RSOS
\maketitle
\else
\par\vspace{1ex}\noindent{\bf Keywords.} \keywordlist.
\fi\fi

\section{Introduction}
The standard approach to the Cayley-Dickson construction~\cite{Harvey} at step $n$ is to define a multiplication operation for algebra $\AA{n} \cong \AA{n-1}\otimes \Z_2$ by
\begin{equation}\label{eqn:cd} (a, b)(c, d) = (ac -\epsilon d^{*}b, da+bc^{*})\ec \end{equation}
where the conjugate is $(a,b)^{*} = (a^{*}, -b)$ and $\epsilon = 1$ for the normed composition algebras analysed by Hurwitz~\cite{Dray}. The products $ac$,  $d^{*}b$, $da$ and $bc^{*}$ are defined according to the construction for $\AA{n-1} \cong \AA{n-2}\oplus \AA{n-2}$, which in turn is eventually determined by elements in $\AA0=\R$.

Traditionally, the multiplication operator is introduced using anticommuting imaginaries $i, j, l$ for the complex numbers, quaternions and octonions, respectively, as
\begin{equation*}\begin{split}
    \AA1 &= \R\oplus i\R = \C\ec \\
    \AA2 &= \C\oplus \C j = \HB\text{ and} \\
    \AA3 &= \HB\oplus l\HB = \OB\es
\end{split}\end{equation*}
This paper introduces a graded decomposition of anticommuting elements, $\o{n}$, at each level, $n$, as 
\[ \AA{n}=\AA{n-1}\oplus \AA{n-1}\o{n} \]
so that $\o1=i$, $\o2=j$ and $\o3=-l$. The ungraded notation of using $k=ij$ is then represented in graded form as $ij=\o{12}$. All elements of $\AA{n}$ are generated using the binomial expansion, which provides an ordering of these elements. This is not important for the normed division algebras but the ordering is found to enable power-associativity, which only comes into play for $\AA{n}$ with $n\ge4$ where $\AA{4} = \SB$, the sedenions.

Deriving an ordering using the graded notation finds that power-associative algebras have three classes of associativity, which leads to four types of non-associativity. Octonions only have one symmetric type of non-associativity that makes the three Moufang loop identities identical. Sedenions introduce two more non-associative types, which separates two of the Moufang identities and allows for zero divisors. The trigintaduonions, $\AA5$, include the final non-associative type, identified by the remaining Moufang identity giving unique results, which adds more complicated zero divisors. Using the new non-associative types provides a simpler categorisation than the Moufang identities. Analysis of the new non-associative types uncovers cyclic forms and modes that categorise zero divisors and applied to the 84 known sedenion zero divisors~\cite{Cawagas1}, reduces them to seven distinct triples of basis elements with three cycles and four modes each. The zero divisors considered here are pairs of the sum of two pure basis terms and linear combinations of these are not considered. The trigintaduonions have 1,260 zero divisors, which are derived from 147 distinct triples via cycles and modes. The new non-associative types identify octonion and three different types of power-associative subalgebras of $\AA{n}$, $n\ge4$, called quasi-octonion algebras in \cite{Wilmot}. Each of these algebras contain twelve zero divisors and since they are embedded in $\AA{n}$, $n\ge4$, in multiples of seven, then it is found that zero divisors in $\AA{n}$, $n>3$ occur in multiples of 84. The cardinality formula for zero divisors for all $\AA{n}$ is derived.

Using loop analysis, Cawagas~\cite{Cawagas1} has previously identified the first power-associative subalgebra by identifying that it did not satisfy one of the three Moufang identities, unlike octonions. He called this quasi-octonions and provided an enumeration in~\cite{Cawagas2}. This paper takes a different approach to loop analysis to prove that most of the sedenion triple products satisfy the three Moufang identities in different ways. Previous work \cite{Wilmot} has shown that the quasi-octonions are part of a series of six algebras, also called quasi-octonions and designated $\PB{k}$, $k\in\{4,8,10,12,14,16\}$ where $k$ is the number of distinct non-associative triple products. These algebras were shown to break the symmetry of octonions in a way consistent with exceptional Lie algebra $\GT$~\cite{Wilmot}. Cawagas found seven copies of $\PB4$ embedded in sedenions as well as eight copies of octonions. Another two of these power-associative algebras, $\PB{12}$ and $\PB{14}$, are introduced in the next two levels of Cayley-Dickson algebras, $\AA5$ and $\AA6$ and only these three $\PB{k}$, and octonions, are embedded in $\AA{n}$, $n\ge6$. Using ${\rm gen}(\PB{k})$ to be a 3-form representation that generates the ultracomplex quasi-octonion, $\PB{k}$, then
\[ \begin{split}
  {\rm gen}(\PB4)   &=\e{123}+\e{145}+\e{167}+\e{246}-\e{257}+\e{347}+\e{356}\ec\\
  {\rm gen}(\PB{12}) &=\e{123}+\e{145}+\e{167}+\e{246}+\e{257}-\e{347}+\e{356}\text{ and}\\
  {\rm gen}(\PB{14}) &=\e{123}+\e{145}+\e{167}+\e{246}+\e{257}+\e{347}+\e{356}\es
\end{split} \]

The octonions have been found to apply to areas in particle physics~\cite{Furey1,Furey2} and extended to sedenions more recently~\cite{Gresnigt1,Gresnigt2,Tang}. It may not be necessary to go to higher dimensions to achieve new representations in particle physics. All Cayley-Dickson algebras of dimension higher than octonions consist of multiple copies of octonions and various quasi-octonions so it is possible that just octonions and the quasi-octonion algebras are all that is necessary.

Section~2 starts by building the Cayley-Dickson construction as two equivalent pyramids using graded and binary notations. Swapping between the two allows the usual composition properties to be proved for all Cayley-Dickson algebras showing that the graded notation is indeed equivalent to the usual doubling of algebras~\cite{Conway}. Section~3 derives the four non-associativity types that simplify the Moufang identities and Malcev algebras and derives the eight non-associative 3-cycles and 3-cycle structure in the Cycle theorem that identifies the subalgebras of $\AA{n}$ for $n\ge3$. The 15 octonion and $\PB4$ representations embedded in sedenions are identified using the graded notation.

Section~4 defines zero divisors and modes in the Mode theorem and shows the relationship to the eight 3-cycles. It proceeds to document the seven and 147 prime zero divisors of $\AA4$ and $\AA5$, respectively, and applies the Cycles and Modes theorems to uncover the rest. Section~5 analyses split Cayley-Dickson algebras, which have the same non-associative structure but different zero divisors. These are provided for the split octonions and both split sedenions.

\section{Graded Cayley-Dickson Construction}
\begin{definition} The basis of the Cayley-Dickson algebra at level $n$, matching the number of generators, has $2^n$ graded elements consisting of {\it grades} of size $\binom{n}k$ called a $k$-grade. Basis elements are represented as ${\rm\bf o_{\alpha}}$ where $\alpha$ is a set of digits from 
$$\MA{n}=\{\emptyset,(1),(12),\ldots,(12\cdots n)\}\ec$$
generated by the recurrence relation $\MA{n}=\MA{n-1} +\MA{n-1}(n)$, equivalent to multiplying the basis $\o{\alpha}$, $\alpha\in\MA{n-1}$ by the generator of $\MA{n}$, $\o{n}$, from the right. Multiplying any basis elements corresponds to a set XOR operation, $\o\alpha\o\beta=\pm\o\delta$, where $\delta=\alpha\veebar\beta$, $\alpha,\beta\in\MA{n}$. The sign is governed by the Cayley-Dickson doubling process and the recurrence relation generates the binomial expansion with the number of terms given by $\sum_{k=0}^n \binom{n}k=2^n$. This is enumerated by Pascal's triangle with recurrence relation, $\binom{n}{k}=\binom{n-1}{k-1}+\binom{n-1}{k}$, where $k\in\N_1^n$, using the notation of Porteous~\cite{Porteous} whereby $\N_1^n$ is the set of natural numbers from 1 up to and including $n$.
\end{definition}

But the Pascal recurrence relation has Grassmannian ordering which is not appropriate for the binary expansion of Cayley-Dickson algebras. This has a pyramid structure whereby the right-hand side repeats the structure on of the left side but multiplied by $\o{n}$. This symmetry enables proof by induction to be visualised, especially when mapped to binary number form. The binary theorem graded pyramid to level 3, starting with the empty set denoting the unit graded element, $\o{\emptyset}=1$, is
\begin{center}
\begin{tabular}{c}
  $1$  \\
  $\{1,\o1\}$  \\
  $\{1,\o1,\o2,\o{12}\}$  \\
  $\{1,\o1,\o2,\o{12},\o3,\o{13},\o{23},\o{123}\}$\es
\end{tabular}
\end{center}
\begin{lemma}\label{lem:1}
Any two distinct graded elements multiplied together generates another distinct graded element.
\end{lemma}
\begin{proof} The XOR operation ($\veebar$) involved in the product corresponds to union ($\cup$) minus intersection ($\cap$) of two sets and has the cancellation property that if $\delta=\alpha\veebar\beta$ then $\delta\veebar\alpha=\beta$ and $\delta\veebar\beta=\alpha$. The XOR operation will cancel any two subsets that are the same, reducing the result to an independent set. This can not be the null set since the elements are distinct so the sets are disjoint.
\end{proof}
\begin{definition} A $k$-grade contains graded elements $\o\alpha$ where $\alpha$ has $k$ indices called the length of $\alpha$. The $k=1$ case is called a {\it generator} because all graded elements can be generated from these by XOR addition. At level $n$ these are $\o1,\o2,\cdots,\o{n}$.
\end{definition}
\begin{definition} A pure element excludes the unit basis element, 1. A {\it blade} is a pure, non-zero, scaled basis element and a $k$-blade is a member of the $k$-grade. Note that blades may have positive or negative squares and for the remainder this is assumed to be $\pm1$, for blades. A blade is positive or negative for scale +1 or -1, respectively. For positive scale, blades have order given by the binomial expansion with $\o2<\o{12}<\o3$, for example. The negative signature blades are smaller than positive blades.
\end{definition}
\begin{definition} A set of blades will be called {\it grade independent} if the corresponding set of indices are linearly independent under the XOR operation.
\end{definition}
The blade $\o\alpha$ represents a pyramid base of $2^n$ numbers one of them being the unit, $\pm1$, and the rest zero. The scalar basis has a unit in the first position and all other positions represent pure basis elements. The usual notation $(a,b)$ translates to $\o\alpha+\o{\beta\cup\{n\}}$ where $a=\o\alpha$ is all zeros in the right-hand half of the pyramid and $b=\o\beta$ and $(\beta\cup\{n\})$ is all zeros in the left-hand half of the $2^n$ numbers.
\begin{definition} The {\it pyramid element} $\va{}$, at level $n$, specifies the $n^{\rm th}$ row of a binary pyramid with coefficients $\va{\beta}$, $\beta\in\MB{n}$, where $\MB{n} = \MB{n-1}(0) \cup \MB{n-1}(1)$ so that $\MB{n}$ is all binary numbers with $n$ digits from $\N_0^1$ and set multiplication extends the binary numbers in $\MB{n-1}$. Hence $\MB{n}$ covers all $2^n$ combinations of $n$ digits in $\MA{n}$ and this is an equivalence relationship between $\o\alpha$ and $\va{\beta}$, $\alpha\in\MA{n}$, $\beta\in\MB{n}$, such that if $i\in\alpha$ then the $i$\textsuperscript{th} position or bit of $\beta$ is 1. Multiplication of $\va{\alpha}$ and $\vb{\beta}$, $\alpha,\beta\in\MB{n}$ involves the XOR operation $\alpha\veebar\beta$, which is bitwise addition modulo 2. At level 1, $a_0=(\pm1,0)=\pm1$ and $a_1=(0,\pm1)=\pm\o1$. To level three this is
\begin{align}
 \va{\alpha} &= (a_0,a_1) \nonumber \\
  &= ((a_{00},a_{10}),(a_{01},a_{11})) \label{eqn:ap2}\\
  &= (((a_{000},a_{100}),(a_{010},a_{110})),((a_{001},a_{101}),(a_{011},a_{111})))\es\nonumber
\end{align}
\end{definition}
The Cayley-Dickson process in pyramid form provides the same ordering as the graded elements with the last index corresponding to $\pm\o{n}$. Some pyramid coefficients at level 2 are $\beta=ij$, $i,j\in\N_0^1$ which corresponds to signed $o_\beta = \o1\ec\:\o2\ec\:\o{12}$ for $a_{ij} = a_{10}\ec\:a_{01}\ec\:a_{11}$ in (\ref{eqn:ap2}), respectively. The four extreme blades at level 3 are
\[ \begin{split}
  \pm\o1 &= (((0,\pm1),(0,0)),((0,0),(0,0))) \text{ and}  \\
  \pm\o{123} &= (((0,0),(0,0)),((0,0),(0,\pm1)))\es
\end{split} \]

\begin{definition}
    The graded version of this Cayley-Dickson construction is more explicit and could be used as a specification for a computer program. For $\alpha,\beta,\delta\in\MB{n}$ then the product $a_\alpha=b_\beta d_\delta$ is the logical AND ($\wedge$) of the bottom level of the pyramid,
    \begin{equation}\label{eqn:cd0}
      a_\alpha = \bigwedge_{i=n,-1}^1 (a_\salpha0i+a_\salpha1i)\ec
    \end{equation}
    where $\alpha_i$ are the bits before the $i$\textsuperscript{th} bit, $\coindex{\alpha_i}$ are the bits after the $i$\textsuperscript{th} position in $\alpha$ and conjugation acts on $\alpha_i$ only. The logical AND product keeps the $\delta$ bits the same at each level and accumulates the sign of each term. Since the $i$\textsuperscript{th} bit is either zero or one then one of the terms in (\ref{eqn:cd}) is always zero and both terms are generated from the Cayley-Dickson XOR construction with the $i$\textsuperscript{th} bit of $\alpha$ and $\beta$ both being the same or both different, respectively, as
    \begin{subequations}
    \label{eqn:cd1}
    \begin{align}
   a_\salpha0i &=  b_\sbeta0i d_\sdelta0i -\epsilon_{i} d^*_\sdelta1i b_\sbeta1i\quad\text{or}\label{eqn:cd1a}\\
   a_\salpha1i &=  d_\sdelta1i b_\sbeta0i +b_\sbeta1i d^*_\sdelta0i\ec \label{eqn:cd1b}
   \end{align}
   \end{subequations}
   where $\beta_i$, $\delta_i$ are the bits yet to be processed and $\coindex{\beta_i}$, $\coindex{\delta_i}$ are the bits already processed in (\ref{eqn:cd0}) and the conjugations only apply to $\alpha_i$ indices.
The (\ref{eqn:cd1a}) terms are the parallel terms and (\ref{eqn:cd1b}) the cross terms. This bit notation is seen as identical to (\ref{eqn:cd}) and it is convenient to sometimes mix the two notations with the $\coindex{\alpha_{i}}$ indices being omitted if empty.
\end{definition}

\begin{definition}
The definition of {\bf conjugation} is $a_\alpha^* = \bigwedge_{i=n,-1}^1 (a_\salpha0i\pm a_\salpha1i)$ with the minus sign selected if $\coindex{\alpha_i}$ is all zeros.
\end{definition}

\begin{lemma} \label{lem:3}
For $\delta\in\MA{n-1}$ not all zeros then $\o\delta\o{n} = \o\beta$, $\beta=\delta\cup\{n\}$. This defines $\o{i_1i_2\cdots i_n}=\o{i_1}\o{i_2}\cdots\o{i_n}$ as positive for $i_1<i_2<\cdots<i_n$.
\end{lemma}
\begin{proof}
  The standard Cayley-Dickson construction is simply 
  \ifnum1=\AACA\hfil\break\fi
  $(\o\alpha,0)(0,1)=(0,\o\alpha)=+\o\beta$ where $\beta=\alpha\cup\{n\}$, which is the same as (\ref{eqn:cd1a}).
\end{proof}

\begin{lemma}\label{lem:2}
For blade $\o{ i_1i_2\ldots i_k}$ then $\o{i_1i_2\ldots i_k}^2=-\epsilon_{i_1}\epsilon_{i_2}\ldots \epsilon_{k}$, for $i_1i_2\ldots i_k, k\in\N_1^n$.
\end{lemma}
\begin{proof}
  For blade $\o\alpha$ with $\alpha\in\MA{n-1}$ then at level $n$ the square of $\o\alpha\o{n}$ from (\ref{eqn:cd1a}) is
  $(0,\o\alpha)(0,\o\alpha)=(-\epsilon_n \o\alpha^*\o\alpha, 0)$, as defined by Lemma~\ref{lem:3}.
  The next iteration of (\ref{eqn:cd0}), if $n-1\in\alpha$ then the result at level $n-1$ is $-\epsilon_{n-1}\epsilon_n \o\alpha^*\o\alpha$ because $\o\alpha^*=-\o\alpha$. If $n-1\notin\alpha$ then the result uses the first term of (\ref{eqn:cd1a}), which repeats the same expression with $a_\alpha^* a_\alpha$ to be evaluated at $i=n-2$. It is easy to see that an inductive argument could successfully be applied to generate the statement of the lemma. This is a well known result proved more generally by Schafer~\cite{Schafer2} and presented here to show the equivalence of the graded notation.
\end{proof}

\begin{theorem}[Cayley-Dickson algebras]\label{thm:cd}
The Cayley-Dickson construction to any level, $n>1$, generates an algebra with pure basis elements, $x$, $y$ satisfying $xy=-yx$ if $x\ne y$ and, if $\epsilon_i=1$, $\forall i\in\N_1^n$, then $xx = -1$. The unit, $+1$, has no effect under multiplication.
\end{theorem}
\begin{proof} 
  Equations (\ref{eqn:cd1}) with reversed terms provide
\begin{align*}
  (b_\sbeta{0}i, &b_\sbeta{1}i) (a_\salpha0i, a_\salpha1i)    \\
    &= \Big(b_\sbeta0i a_\salpha0i -\epsilon_{i} a^*_\salpha1i b_\sbeta1i,
        \: a_\salpha1i b_\sbeta0i +b_\sbeta1i a^*_\salpha0i \Big) \\
    &= -\Big(a_\salpha0i b_\sbeta0i -\epsilon_{i} b^*_\sbeta1i a_\salpha1i,
        \: b_\sbeta1i a_\salpha0i +a_\salpha1i b^*_\sbeta0i \Big) \\
    &= -(a_\salpha0i, a_\salpha1i)(b_\sbeta{0}i, b_\sbeta{1}i)\es
  \end{align*}
  This is a result proved by Schafer~\cite{Schafer2} and this proof is used to show the equivalence of the graded notation. The result for $x^2=-1$ is from Lemma~\ref{lem:2} with $\epsilon_i=1$, $\forall i\in\N_1^n$, and the unit statement is obvious.
\end{proof}

\begin{corollary}
The first three levels of blades with $\epsilon_i=1\;\forall i\in\N_1^3$ are isomorphic to complex numbers, $\C$, quaternions, $\HB$, and octonions, $\OB$, respectively.
\end{corollary}
\begin{proof}
Theorem~\ref{thm:cd} showed that blades generate Cayley-Dickson algebras under the product rules (\ref{eqn:cd1}). It is obvious that the first generator, $\o1$, at level 1 is isomorphic to the imaginary of complex numbers. At level two, the usual representation of quaternions also have the cyclic product rules, $ij=k$, $jk=i$ and $ki=j$. Starting by
\ifnum1=\RSOS
stating the cyclic product rules for quaternions, which are easily proved and follow the usual doubling construction,
\else
calculating the cyclic rules in full
\fi
\begin{equation} \label{eqn:qc3}
\ifnum1=\RSOS
 \o1\o2=\o{12}\ec\quad\o2\o{12}=\o1\quad\text{and}\quad\o{12}\o1=\o2\es
\else
 \begin{split}
  \o1\o{2}  &= ((0,1),(0,0))((0,0),(1,0)) = ((0,0),(0,1)(1,0)) \\
            &= ((0,0),(0,1^*)) = \o{12}\ec \\
  \o2\o{12} &= ((0,0),(1,0))((0,0),(0,1)) = (-(0,1)^*(1,0),(0,0)) \\
            &= ((0,1),(0,0)) = \o1 \ \text{ and}\\
  \o{12}\o1 &= ((0,0),(0,1))((0,1),(0,0)) = ((0,0),(0,1)(0,1)^*) \\
            &= ((0,0),-(-1,0)) = \o2\es
 \end{split}
\fi
\end{equation}
Thus level 2 is isomorphic to $\HB$. At level 3 there are seven blades, $\{\o1, \o2, \o{12}, \o3$, $\o{23}, \o{13}, \o{123}\}$, which matches the number of octonion basis elements and it is obvious that the quaternions are embedded in the first three pure terms, which is a feature of $\OB$. Applying Artin's Theorem~\cite{Harvey}, whereby any pair of basis elements generates an associative subalgebra, the following rings can be identified,
\begin{equation}\begin{split}\label{eqn:qa7}
     &(\o1, \o2, \o{12}), (\o2, \o3, \o{23}), (\o1, \o3, \o{13}), (\o1, \o{23}, \o{123}),\\
     &(\o2, \o{13}, \o{123}), (\o{12}, \o3, \o{123}), (\o{12}, \o{13}, \o{23})\es
\end{split}\end{equation}
These are quaternion apart from the fourth and seventh terms, which are anti-quaternion, and swapping the order of two terms makes these quaternion rings. A calibration of seven quaternions specifies an algebra isomorphic to octonions because they satisfy the rules of Theorem~\ref{thm:cd}.
\end{proof}

The alternative loop properties of Cayley-Dickson multiplication that are essential to simplify non-associative expressions are now proved~\cite{Bruck1}. The notation of left-evaluation or right-expansion is introduced to make it easier to perform all proofs. This is the natural processing order for computer interfaces so integrates with existing parsers naturally.
\begin{lemma} \label{lem:4}
For $x$ and $y$ distinct blades with $\epsilon_i=1$, $i\in\N_1^n$, then $yxx = -y$ and $x(xy) = -y$, using left-evaluation.
\end{lemma}
\begin{proof}
Schafer~\cite{Schafer2} proves the first statement in his Lemma 4 and uses the flexibility property~\cite{Schafer1}, $xyx=x(yx)$, to prove the second. This property is derived from the anticommutivity of Theorem~\ref{thm:cd}, for distinct blades $a,b,c$ then $abc=c(ba)$. Schafer's proof is very brief and it is necessary that the same inductive approach works for the graded approach. Here the former statement, $yxx=-y$ is proved following the same inductive argument as Schafer, and the flexibility property shows the two statements are equivalent. It works for $\HB$ which needs to be checked for the six combinations of two from three elements using the cyclic product rules (\ref{eqn:qc3}) and the $xx=-1$ property from Theorem~\ref{thm:cd},
\begin{align*}
  \o1\o2\o2       &= \o{12}\o2=-\o2\o{12}=-\o1\ec \\
  \o1\o{12}\o{12} &= -\o2\o{12} = \o{12}\o2 = -\o1\ec \\
  \o2\o1\o1       &= -\o1\o2\o1=-\o{12}\o1=-\o2\ec \\
  \o2\o{12}\o{12} &= -\o{12}\o2\o{12} = \o1\o{12} = -\o2\ec \\
  \o{12}\o1\o1    &= -\o1\o{12}\o1 = \o2\o1 = -\o{12} \text{ and} \\
  \o{12}\o2\o2    &= -\o1\o2 = -\o{12}\es
\end{align*}
Induction means proving expressions involving $\o{n}$ assuming $\o\alpha\o\beta\o\beta=-\o\alpha$, for $n\notin\alpha$ or $\beta$. After the first step of (\ref{eqn:cd1}), the path through the pyramid on the right-hand side with $\o{n}$ is the same as the path on the left-hand side without $\o{n}$ so the inductive argument invokes the symmetry of the pyramid. Assuming $\va{\alpha}\vb{\beta}\vb{\beta}=-\va{\alpha}$, $\alpha, \beta\in\MB{n-1}$ not all zeros, then only the following conditions need to be proved,
\begin{align*}
 \va{\alpha1}\vb{\beta}\vb{\beta} &= (0,a_{\alpha})(b_{\beta},0)(b_{\beta},0) \\
   &= (0,a_{\alpha}b_{\beta}^*)(b_{\beta},0) \\
   &= (0,a_{\alpha}b_{\beta}^*b_{\beta}^*) \\
   &= -\va{\alpha1}\ec
\end{align*}
\begin{align*}
 \va{\alpha}\vb{\beta1}\vb{\beta1} &= (a_{\alpha},0)(0,b_{\beta})(0,b_{\beta}) \\
   &= (0,b_{\beta}a_{\alpha})(0,b_{\beta}) \\
   &= (-b_{\beta}^*(b_{\beta}a_{\alpha}),0) \\
   &= (b_{\beta}(b_{\beta}a_{\alpha}),0) \\
   &= -\va{\alpha}\quad \text{ and }
\end{align*}
\begin{align*}
 \va{\alpha1}\vb{\beta1}\vb{\beta1} &= (0,a_{\alpha})(0,b_{\beta})(0,b_{\beta}) \\
   &= (-b_{\beta}^*a_{\alpha},0)(0,b_{\beta}) \\
   &= (0,-b_{\beta}(b_{\beta}^*a_{\alpha})) \\
   &= (0,b_{\beta}(b_{\beta}a_{\alpha})) \\
   &= (0,-a_{\alpha}) \\
   &= -\va{\alpha1}\es
\end{align*}
This leaves $y\o{n}\o{n}$, which follows from the last two equations with $\beta$ all zeros or simply 
$(a,1)(0,1)=(-\epsilon_n1,a)=-a$. Hence $\va{\alpha}\vb{\beta}\vb{\beta}=-\va{\alpha}$ for all $\alpha,\beta\in\MB{n}$. Note that while $xyzz = -xy$, it is not necessarily true that $x(xyz)=-yz$ for distinct blades $x, y, z$. For example, $\o1(\o1\o2\o3)=\o{23}$.
\end{proof}

\begin{definition} The usual {\it associator} is defined, using left-evaluation, as
\[ [a, b, c] = (ab)c - a(bc) = abc -a(bc)\es\]
This is designated the {\it ordered associator} for $a<b<c$ and other orderings are called {\it unordered associators}.
\end{definition}
A triple of elements from $\HB$ under multiplication generates a scalar and the associator is zero for quaternion and anti-quaternion rings, $[b,c,bc]=bc(bc)-b(c(bc))=(bc)^2+b(bcc)=0$. Octonions have $\binom73 = 35$ ordered triples with seven being rings, as specified by (\ref{eqn:qa7}). The remaining 28 triples have non-zero associators. For example, Lemma~\ref{lem:2} has $\o1\o2\o3 = \o{123}$ whereas expanding the pyramid for $\o1(\o2\o3)$ gives the term $-a_{100}b_{011}$ in the $\o{123}$ position so that $\o1(\o2\o3) = -\o{123}$, hence the associator is non-zero and this triple is non-associative, 
\ifnum1=\AACA\break\fi $[\o1,\o2,\o3] = 2\o{123}$. The graded notation allows power-associativity to be analysed when the $\o{n}$, $n\ge4$, generators are introduced.

The names quaternion and octonion seem to be derived from book binding techniques involving a quire of four or eight folios of paper, respectively. Traditionally, an ancient book or codex would use separate sheets or bifolio with a single fold stacked and bound to the binding or stub. This was superseded with folding a single large sheet several times to make the quire. The type of fold starts with a bifolio, which is double the number of folios. Two bifolia have quarto fold type, four has Octavo and eight has extodecimo type, which covers more than eight folios~\cite{Deroche}. The power-associative algebras are such a vast departure from hypercomplex numbers that it is necessary to distinguish this departure. For historical reasons the prefix exto has been transformed to the modern term ultra and used to designate all power-associative Cayley-Dickson algebras, including the quasi-Cayley-Dickson algebras. Then, the complex numbers are followed by the hypercomplex numbers, $\HB$ and $\OB$, and the following power-associative algebras would be known as the ultracomplex numbers.
\begin{definition}
The Cayley-Dickson series for power-associative algebras after $\AA3$, with $\epsilon_i=1$ for $i\in\N_1^n$, are designated $\UB{m}$, $m=n-3$, $m>0$, and referred to as {\it standard ultracomplex algebras}. Thus the normed algebras, $\R$, $\C$, $\HB$ and $\OB$, are extended to the $\UB{m}$ series starting with sedenions, $\UB1$. The $\UB{m}$ algebras embed quasi-octonions and reference~\cite{Wilmot} found six quasi-octonions, $\PB{k}$ for $k\in\{4,8,10,12,14,16\}$. This paper finds only three occur in the power-associative Cayley-Dickson algebras and these are denoted $\UB0$. This set of quasi-octonions, along with octonions, are the basis of ultracomplex algebras, and there are quasi-Cayley-Dickson algebras at all levels, $\UB{m}$, $m\ge0$, but only a subset are included in $\UB{m+1}$. It is this subset of all quasi-Cayley-Dickson algebras that are named ultracomplex that embed various combinations of $\UB{0}$.
\end{definition}

\section{Structure Theorems}
To concentrate on non-associative triple products the term triad is used to provide distinct triples that avoid scalars and repeated elements with trivial associativity. The notation of left-evaluation has simplified the product rules to only right-hand multiplication, the first two terms of any expression can be swapped if the sign is changed and the last two basis elements will cancel if equal, apart from sign. Since evaluation is performed from the left then any pair of successive common basis elements inside an expression will cancel.

\begin{definition} There are $N^3$ scaled triple products $bcd$ and $\binom{N}3$ naturally ordered triples with blades $b<c<d$ where $N=2^n-1$ and $n$ is the number of generators. The ordered triples are the same as removing all pairs and triples of the same element from the $N^3$ set and dividing by the number of permutations, $6$. A {\it triad} is defined as ordered triples containing distinct blades whereby $b<c<d$ so that $a=bcd$ can be defined for triads and $a$ is distinct from $b, c, d$, which follows from Lemma~\ref{lem:1}. Of course, there are 5 unordered triples for each triad. The product $a$ will be scalar unless the triad is grade independent.
\end{definition}

\begin{theorem}[Associativity theorem]\label{thm:as}
The $24$ permutations of triad products generate three classes of unordered associativity, called Classes~1, 2 and 3, with Class~1 being dependent on either Class~2 or~3. These are shown in \tar{tat} where equivalence of associators means both are zero or both are non-zero, noting that only 12 permutations for ordered $b, c, d$ and $a=bcd$ are required because $[a,b,c]=-[c,b,a]$.
\tab{ht}{tat}{Triad Unordered Associativity Classes}
\begin{tabular}{|c|c|} \hline
{\bf Equivalent Associators} & {\bf Class} \\\hline
$[b, a, c]\approx[b, d, c]\approx[a, b, d]\approx[a, c, d]$ &Class~1\\\hline
$[a, b, c]\approx[b, c, d]\approx[b, a, d]\approx[a, d, c]$ &Class~2\\\hline
$[a, c, b]\approx[c, b, d]\approx[a, d, b]\approx[c, a, d]$ &Class~3\\\hline
\end{tabular}\tae
\end{theorem}
\begin{proof} Excluding complete associativity for triads where all associativity classes are zero means that $a(bc)=-(bc)a=cba$ if $a$ is not a scaling of $bc$ and similarly for $abc=(ab)c=c(ba)$. Hence $[a,b,c] = abc -a(bc) = c(ba) -(cb)a = -[c,b,a]$. This is also true for associative triple products, where $a$ is a scalar, because both sides are zero. Now, if $b,c,d$ are not associative, the relation $a=bcd$ means $ad^{-1} = bc$ or $d = -a(bc)^{-1} = -bca/b^2/c^2$, noting that $(bc)^2 = -b^2c^2$. First, Class~2 and Class~3 associativity are analysed.

Assuming Class~2 associativity, $[a,b,c] = abc -a(bc) = 0$ or $[c,b,a] = cba -c(ba) = 0$, then sub\-stituting $d$ into the other three associators in the second row in \tar{tat}, scaled by $(bc)^2$, gives
\begin{align*}
     {[}b,c,d](bc)^2 &= bc(bca) -b(c(bca)) = -(bca)(bc) -(bca)cb \\
        &= -(bc)a(bc) + c(ba)cb = a(bc)(bc) -baccb \\
        &= a(bc)^2 +abbc^2 = -ab^2c^2 +ab^2c^2 = 0\ec \\
     {[}b,a,d](bc)^2 &= ba(bca) -b(a(bca)) = -ba(c(ba)) +b(bc)a^2 \\
          &= c(ba)(ba) +cbba^2 = 0\text{\quad and}\\
     [a,d,c](bc)^2 &= a(bca)c -a(bcac) = -(bca)ac +a(c(ba)c) \\
         &= -bcca^2 +(bacc)a = 0\es
\end{align*}

Assuming Class~3 associativity, $[a,c,b] = acb -a(cb) = 0$ or $[b,c,a] = bca -b(ca) = 0$, then sub\-stituting $d$ into the other three associators in the third row of \tar{tat}, scaled by $(bc)^2$, gives
\begin{align*}
     [c,b,d](bc)^2 &= cb(bca) -c(b(bca)) = -b(ca)(cb) -bcabc = cab(cb) -b(ca)bc \\
                   &= -a(cb)(cb) +(ca)bbc = a -cac = 0\\
     {[}a,d,b](bc)^2 &= a(bca)b -a(bcab) = -bcaab +bcaba = bcb +b(ca)ba = 0\text{\; and}\\
     {[}c,a,d](bc)^2 &= ca(bca) -c(a(bca)) = -bca(ca) -bcaac = -b(ca)(ca) +bcc = 0\es
\end{align*}

Assuming Class~1 associativity, $[b,a,c] = bac -b(ac) = 0$ or $[c,a,b] = cab -c(ab) = 0$, then substituting $d$ into the other three associators in the first row of \tar{tat}, scaled by $(bc)^2$, gives
\begin{align*}
     {[}b,d,c](bc)^2 &= b(bca)c -b(bcac) = -bcabc +bcacb \ec \\
     [a,b,d](bc)^2 &= ab(bca) -a(b(bca)) = (bca)(ba) -bcaba = -cba(ba) +cbaba\text{ and} \\
     {[}a,c,d](bc)^2 &= ac(bca) -a(c(bca)) = -bca(ac) -bcaca\es
\end{align*}
 If Class~2 associativity is also assumed then $cba=c(ba)$ and, also again using $[c,a,b]=0$, the three equations become
\begin{align*}
     -bcabc +bcacb   &= c(ba)bc -c(ba)cb = -b(ac)bc +baccb = acbbc +abb = 0\ec \\
     -cba(ba) +cbaba &= c(ba)(ba) -c(ba)ba = -c +bac(ab) \\
          &= -c +b(ac)ba = -c +cabba = 0\text{ and} \\
     -bca(ac) -bcaca &= c(ba)(ac) +c(ba)ca = -bac(ac) -bacca \\
          &= -b(ac)(ac) -b = 0\es
\end{align*}
If Class~3 associativity is instead assumed with Class~1 then $bca=b(ca)$ or $cba=-b(ca)$ and, also again using $[c,a,b]=0$, the three equations become
\begin{align*}
     -bcabc +bcacb   &= -b(ca)bc +b(ca)cb = -acbbc -cabcb \\
         &= -a -c(ab)cb = -a +abccb = 0\ec \\
     -cba(ba) +cbaba &= b(ca)(ba) -b(ca)ba = -cab(ba) +cabba \\
         &= c(ab)(ab) +c = 0\text{ and} \\
     -bca(ac) -bcaca &= b(ca)(ca) -b(ca)ca = -b +cabca \\
         &= -b +c(ab)ca = -b +bacca = 0\es
\end{align*}
Hence Classes~2 and 3 associativity are non-associatively equivalent in their own right but Class~1 needs either Class~2 or Class~3 for equivalence. Equivalence classes do not need to be independent and the non-associativity theorem below explains why this connection occurs.
\end{proof}

\begin{definition}
The {\it triple associator}, $T(b,c,d)$, is introduced to cover associativity Class\-es~1, 2 and~3, which allows for an unambiguous definition of associativity,
\[ T(b,c,d) = [b, d, c] -[d, c, b] +[c, b, d]\es \]
\end{definition}
\begin{definition}A triad, $b,c,d$, is {\it associative} if 
\ifnum1=\ELS\break\fi
$T(b,c,d)=0$ and {\it non-assoc\-iative} if \ifnum1=\ARXIV\break\fi $T(b,c,d)\ne0$. This may seem like over-kill but will soon be seen to be a necessary condition.
\end{definition}

Artin's theorem~\cite{Harvey} is stated only for associative subalgebras of octonions but can be gen\-er\-alised and then extended to find all subalgebras of any Cayley-Dickson algebra.
\begin{theorem}[Generalised Artin theorem]\label{thm:at}
Any two distinct blades of $\HB$, $\OB$ or $\UB{m}$, $m>0$, generate an algebra isomorphic or anti-isomorphic to quaternions, which is thus associative. The map $\langle \phi :\:q\rightarrow-q\rangle$, $q\in\HB$ generates anti-quaternions whereby $\phi(q_i)\phi(q_j)=-\phi(q_i q_j)$, $q_i\in\HB$, $i\ne j$.
\end{theorem}
\begin{proof}
From Lemma~\ref{lem:1}, multiplying any two distinct blades generates a third distinct blade and this forms a (non-Abelian) ring by the XOR cancellation property. All such products are closed under multiplication and form a Cayley-Dickson algebra, by Theorem~\ref{thm:cd}, because distinct blades anti-commute and have negative square in this case. These rings show the cyclic quaternion product rules, since for any two blades, $b$ and $c$ then $bc=d$, $cd=c(bc)=-bcc=b$ and $db=-cbb=c$, by Lemma~\ref{lem:4}, and thus the triad $b$, $c$ and $bc$ is isomorphic or anti-isomorphic to $\HB$, depending on whether $bc(bc)=-1$ or $1$, respectively. Swapping $b$ and $c$ for the latter case changes the triple product to $-1$, producing quaternions. The triple associator shows associativity for this triad,
\[ \begin{split} T(b,c,bc) &= b(bc)c -b(bcc) -bc(cb) +bc(cb) +cb(bc) -c(b(bc)) \\
                           &= -bcbc +bb -bc(bc) +c(bcb) \\
                           &= +1 -1 +1 -1 \\
                           &= 0\es \\
\end{split} \]
\end{proof}
\begin{theorem}[Cayley-Dickson algebra generation]\label{thm:ug}
Any power-associative Cayley-Dick\-son algebra, $\UB{m}$, $m>0$, as well as $\OB$, can be generated by $n=m+3$ distinct blades, $\o{\alpha_k}, \forall\:k\in\N_1^n$, such that $k\in\alpha_k$ and the blades are grade independent.
\end{theorem}
\begin{proof}
 For convenience, temporarily allow $\UB0$ to include $\OB$. Triad $b,c,d$ (ie three distinct blades), with $a=bcd$ non-scalar, is non-associative,
\begin{equation}\begin{split}
  T(b,c,d) &= bdc -b(dc) -dcb +d(cb) +cbd -c(bd) \\
           &= bdc +dcb -dcb +bcd -bcd +bdc \\
           &= 2bdc \\
           &= \pm2a\es \label{eqn:jna}
\end{split}\end{equation}
Each row of Pascal's triangle of graded elements contains an even number of each index, excluding level 1, because each index is doubled for the next level of the recurrence relation and there are an even number of existing elements to double with the new index. Hence, excluding the first row, the multiplication of all elements of any row cancels the entire set because they are grade dependent. The $\C$, $\HB$ and $\OB$ cases show that an element can be selected from the row of Pascal's triangle that is independent from all previous rows. At every level this can be done by choosing an element from the right-hand-side of the pyramid that contains $\o{n}$ and whose other indices are linearly independent under the XOR operation of the set chosen from the previous row. Dependency would generate a whole row of Pascal's triangle and cancel the generator at that level. Excluding such cases produces a distinct new graded element.

Hence, any set of basis elements, $\o{\alpha_k}$, $k\in\alpha_k$, that are grade independent in the $n=m+3$ generators, $\o{k}$, will expand the basis to $N_n=2^n-1$ distinct pure basis elements by applying Lemma~\ref{lem:1} repeatedly. By definition all triples of this set are triads, which by (\ref{eqn:jna}) are non-associative. Any pair of elements generates a distinct third element by Lemma~\ref{lem:1}, which is a ring by Theorem~\ref{thm:at} and terms can be swapped in order to make these quaternion. The number of quaternion subalgebras, $H_n$ is given by selecting two of the $N_n$ basis elements, which generates a third and this can be done in six ways, so
\begin{equation} 
   H_n = \frac16 N_n(N_n -1)\es\label{eqn:hn}
\end{equation}
At level 3, combinations of three blades, $b,c,d$, generate seven distinct blades with, $bc, bd, cd, bcd$, by Lemma~\ref{lem:1}. By Theorem~\ref{thm:cd} these satisfy the Cayley-Dickson properties and hence generate $H_3=7$ quaternions, by Theorem~\ref{thm:at}. The corollary to Theorem~\ref{thm:cd} shows a calibration of seven quaternions derived from $\OB$ and this is the reverse situation where $b,c,d$ with $bcd$ non-scalar generate seven quaternions that are isomorphic to $\OB$, if $\HB$ is embedded first. At level 4, $H_{15}=35$ quaternions are generated by the $N_{4}=15$ basis elements as provided by the calibration in \cite{Wilmot}, which extends the octonions calibration thus generating an algebra isomorphic to sedenions, $\SB$.  

Equation (\ref{eqn:hn}) defines the number of quaternions needed to define an associative calibration whereby each pair of elements is included once. A set of $H_n$ distinct pairs can be selected and by Theorem~\ref{thm:at} extended to triads, terms swapped if necessary for quaternions, to generate a calibration that specifies a product table for all pairs of elements. The quaternions can be selected so that $\UB{m-1}$ is embedded as a subalgebra in the first part of the algebra, $m>0$, and this construction satisfies Theorem~\ref{thm:cd} thus generating an algebra isomorphic to $\UB{m}$. The value of $H_n$ for higher levels is demonstrated in \tar{usa} later.
\end{proof}
The generation theorem disguises the complexity of the algebras since it avoids analysing the non-associativity of arbitrary triads. Analysing the asymmetry of non-associativity uncovers the triads related to power-associativity.
\begin{theorem}[Non-associativity theorem]\label{thm:na}
If a triad is not associative then either all associativity classes are non-zero or only one is non-zero, the other two being zero.
\end{theorem}
\begin{proof}
From \tar{tat}, for triad $b,c,d$, Class~1 associativity has $bdc=b(dc)$, Class~2 associativity has $bcd = b(cd)$ and Class~3 associativity has $cbd=c(bd)$. Starting with Class~1 and substituting Class~3 then Class~2 associativity, is
\[ bdc = -c(bd) =  -cbd = bcd = -b(dc)\ec \]
which contradicts Class~1 associativity if Classes~2 and 3 are assumed. Similarly, substituting the reverse of Class~3 then Class~1 into Class~2,
\[ bcd = -d(bc) = -dbc = bdc = b(dc) = -b(cd)\ec \]
which contradicts Class~2 associativity. Finally substituting the reverse of both  Class~2 and  Class~1 into  Class~3,
\[ cbd = -d(cb) = -dcb = cdb = c(db) = -c(bd)\ec \]
which contradicts Class~3 associativity. This proves existence of associativity pairs for non-associativity, which later are used to classify the types of power-associativity. For necessity of non-associativity, since $a=bcd$ then each associator in (\ref{eqn:jna}) is 0 or $\pm2a$ and since the sum is $\pm2a$ then either exactly one or all three associators are non-zero.
Hence the definition of non-associativity as $T(b,c,d)\ne0$ is adequate.
\end{proof}

\begin{definition} 
{\it Non-associativity types} are defined, using \nz\  to mean non-zero, as
\begin{center}
\begin{tabular}{|c|c|c|c|} \hline
\multicolumn{3}{|c|}{\bf Ordered Associativity} & \multirow{2}{*}{\bf Non-Associativity} \\\cline{1-3}
{\bf Class~1} &{\bf Class~2} &{\bf Class~3} &  \\\hline
 \nz  &0   &0   &{\bf Type~A}\\\hline
 0    &\nz &0   &{\bf Type~B}\\\hline
 0    &0   &\nz &{\bf Type~C}\\\hline
 \nz  &\nz &\nz &{\bf Type~X}\\\hline
\end{tabular}
\end{center}
\end{definition}
\begin{definition}
A triple $(b, c, d)$ generates a {\it Malcev algebra} under the operation $x\ast y=\frac12(xy-yx)$ if it satisfies Malcev's identity~\cite{Sagle}, $J(x,y,[x,z])=[J(x,y,z),x]$, where the Jacobian is
\[ J(x,y,z)=[[x,y],z] +[[y,z],x] +[[z,x],y]\es\]
\end{definition}
\begin{lemma}
  The Malcev identity simplifies for grade independent triads to
\begin{equation}\label{eqn:mai} bcdb+bc(db) = 0\es
\end{equation}
\end{lemma}
\begin{proof}
For the triple being a triad that only involves pure, distinct elements, the commutator is $[x,y]=(xy-yx)=2xy$ so that the Jacobian for grade independent triad $(x,y,z)$, becomes 
\[\begin{split}
  J(x,y,z)&=4\big(xyz-z(xy)+yzx-x(yz)+zxy-y(zx)\big)\\
    &= 4[x,y,z] +4[y,z,x] +4[z,x,y]\es
\end{split}\]
This is similar to the triple associator but with all positive signs. For grade independent triads, the Malcev identity simplifies as
\begin{equation}\begin{split}
  \frac14(&J(b,c,bd)-J(b,c,d)b) \\
   &=bc(bd)-bd(bc)+c(bd)b-b(c(bd))+bdbc-c(bdb)\\
   &-bcdb+d(bc)b-cdbb+b(cd)b-dbcb+c(db)b\\
   &=2bc(bd) -bdcb-bdcb -dbbc-dbbc\\
   &-bcdb+bcdb-cd+cdbb+dbcb+dbcb\\
   &=2bc(bd) -2bdcb -2bcdb -2dbcb\\
   &=-2bc(db) -2bcdb\\
   &=0\es
\end{split}\end{equation}
Note that grade independence allows $c(bd)=-bdc$ since $bdc$ is not a scalar.
\end{proof}
Sagle~\cite{Sagle} provides an equivalent form to this equation but assumes alternate-associ\-ativity to derive it. The identification of the Malcev identity is $[bc,d,b] \ne0$, which implies A or X non-associativity if $bc$ is interpreted as a new $c$. Malcev algebras generalise Lie algebras to non-associative algebras and this is analysed later where B and C non-associativity types that satisfy Malcev's identity are also encountered.

\begin{theorem}[Moufang associativity theorem]\label{thm:ma}
A Moufang loop~\cite{Bruck1,Bruck2} is a non-assoc\-iative algebra that satisfies one of three identities for non-assoc\-iative triad $(b,c,d)$. Calling these types 1 to 3, then each type exhibits the non-associativity shown in the following table where the Moufang Identity condition equalling zero means the triad is non-associative.
\begin{center}
\begin{tabular}{|c|c|c|} \hline
{\bf Type} & \text{\bf Moufang Identity} & {\bf Association} \\\hline
{\rm1} & $d(b(dc)) -dbdc = 0$ & \text{\rm B or X} \\
{\rm2} & $b(d(cd)) -bdcd = 0$ & \text{\rm C or X} \\
\multirow{2}{*}{{\rm3}} & \multirow{2}{*}{$db(cd) -d(bc)d = 0$} & \text{\rm B or X if Malcev,} \\
  & & \text{\rm A or C otherwise}\\\hline
\end{tabular}
\end{center}
\end{theorem}
\begin{proof} Moufang type~1 and 2 equations, multiplied by $d$ and assuming $d^2 = -1$, become respectively
\begin{align*}
     d(b(dc))d -dbdcd &= -b(dc)dd +bddcd = dcbdd -bcd \\
             &= -dcb -bcd = -bcd -b(cd) = 0\ec \\ 
     b(d(cd))d -bdcdd &= -d(cd)bd +dbcdd = -cddbd +dbcdd \\
             &= cbd -dbc = cdb +c(db) = 0\es
\end{align*}
Substituting the first equation into Class~2 associativity is $[b,c,d] = bcd-b(cd) = 2bcd$, which is non-zero. Substituting the second into Class~3 associativity is $[c,b,d] = cbd -c(bd) = 2cbd$, which is also non-zero. Hence Moufang loop non-associativity of Types~1 and 2 correspond to B-non-associativity and C-non-associativity, respectively, or X-non-associativity in both cases.

Multiplying Moufang Type~3 by $db$ gives,
\[ db(cd)(db) -d(bc)d(db) =-cd(db)(db) -bc(db) = 0\ec \]
so $cd = bc(db)$. Applying (\ref{eqn:mai}) so that $cd=-bcdb$ and substituting into Class~2 associativity gives,
\[ [b,c,d] = bcd -b(cd) =bcd +b(bcdb) =bcd -bcdbb =bcd +bcd =2bcd \ne 0\es \]
This implies type B or X non-associativity. If Malcev's identity is not satisfied then (\ref{eqn:mai}) takes the opposite sign, $[bc,d,b]=0$, and Class~2 associativity is then zero so that A or C non-associativity holds. This completes the proof and it is noted that it seems possible to use combinations of the Moufang types along with Malcev's identity to distinguish all four types of non-associativity A, B, C and X. This is the case for types B, C and X, but it is found that some A types do not satisfy the Moufang identity.
\end{proof}

Note that some authors~\cite{Kunen,Vojtechovsky} include a fourth Moufang type, $db(cd) = d(bcd)$. But this is the same as the third type since $d(bcd)=-bcdd=d(bc)d$. Also note that some authors follow Bruck~\cite{Bruck1} and assume that the three Moufang identities are equivalent. This is only true for alternate-associativity and the association table of Theorem~\ref{thm:ma} obviously generalises to define power-associativity. We will find that both Moufang identities and the four non-associative types change for various triads in the subloops of Cayley-Dickson algebras but it is easier to classify the pattern using the four new non-associative types.

The associative triads have pairs $(b,c)$, $(b,bc)$ and $(c,bc)$, which form a 3-cycle of pairs. This applies to all pairs and can be extended to triads to provide another level of structure to the ultracomplex algebras. 
\tab{ht}{tss}{Triad Structure}
\begin{tabular}{|c|c|c|c|c|}  \hline 
{\bf Label} &{\bf Associative} &{\bf 3-Triad Cycles} &{\bf Non-Cycles} &{\bf Total} \\\hline
$\C$ &0 &0 &0 &0 \\
$\HB$ &1 &0 &0 &1 \\
$\OB$ &7 &24 &4 &35 \\
$\UB1$ &35 &360 &60 &455 \\
$\UB2$ &155 &3720 &620 &4,495 \\
$\UB3$ &651 &33,480 &5580 &39,711 \\
$\UB4$ &2667 &283,464 &47,244 &333,375 \\
$\UB5$ &10,795 &2,331,720 &388,620 &2,731,135 \\
$\UB6$ &43,435 &18,912,840 &3,152,140 &22,108,415 \\
$\UB7$ &174,251 &152,345,160 &25,390,860 &177,910,271 \\\hline
\end{tabular}\tae

\begin{definition} Naturally ordered basis pairs form {\it 3-cycles}, $(b,c)$, $(b,bc)$ and $(c,bc)$, which cover all ordered pairs. This is because there are $(2^n-1)(2^n-2)/2$ combinations of pairs which have factors one and two numbers less than a power of $2$ so have a factor of $3$ and the 3-cycle is closed. These pairs can be ordered as $b<c<bc$ and the pairs of the ordered 3-cycles are called the first, second and third cycle, respectively. A {\it 3-triad cycle} extends this to the triads $(b,c,d), (b,bc,d), (c,bc,d)$ where $b<c<bc<d$. These are not distinct since triads with $d<bc$ have second and third cycles already included elsewhere. Such cases are called {\it non-cycles} and degenerate cases, $d=bc$, give the associative triads. This divides the triads into three parts, associative, 3-triad cycles and non-cycles.
\end{definition}

\tar{tss} enumerates all triads for Cayley-Dickson algebras for up to $10$ generators. Non-associative triads are split into non-cycles and cycles where the non-cycle ones have $d<bc$. The Associative column is $H_n$ from (\ref{eqn:hn}) while the Cycles column have been generated by code described in the last section. The summation of these triad enumerations is provided in the last column of \tar{tss} that matches $\binom{N}3$, where $N=2^n-1$. 

For quaternions, $\HB$, the associative triad in \tar{tss} is $(\o1,\o2,\o{12})$. For octonions, $\OB$, the $7$ associative triads have been provided by (\ref{eqn:qa7}). The $24$ 3-triad cycles are derived from eight triads starting with the following first cycles,
\begin{equation}\label{eqn:nox}\begin{tabular}{cccc}
$(\o1, \o2, \o3)$, &$(\o1, \o2, \o{13})$, &$(\o1, \o2, \o{23})$, &$(\o1, \o2, \o{123})$, \\
$(\o1, \o3, \o{23})$, &$(\o1, \o3, \o{123})$, &$(\o2, \o3, \o{123})$,
       &$(\o{12}, \o{13}, \o{123})$.\\
\end{tabular}\end{equation}
All these octonion triads have X-non-associativity and it will soon be proved that the cycles must also be of type X. Remember $[\o1,\o2,\o3] = 2\o{123}$ and multiplying terms by distinct basis elements does not change the non-associativity type. To achieve this we need another generator, $\o4$, for example, so that for $b=\o1$, $c=\o2$ and $d=\o{34}$ then Class~2 is associative, $[\o1,\o2,\o{34}] = -\o{1234} +\o{1234} = 0$ and Class~1 is non-associative, $[\o1,\o{34},\o2] = -\o{1234} -\o{1234} \ne 0$.
\tar{aeg} provides the first A and C triads from $\UB1$ and the first B triad from $\UB2$ as examples. For X non-associativity, the first $\UB1$ X case not included in $\OB$ is selected but note that any octonion non-associative triad could have been selected with generator $\o3$ replaced with $\o4$.
\begin{definition} The {\bf prime} triad is the first ordered triad that generates the loop of the 24 3-triads is called the  triad. For the 35  triads given by (\ref{eqn:qa7}) and the 3-triad cycles of (\ref{eqn:nox}), this is the first cycle in (\ref{eqn:qa7}). By  Theorem~\ref{thm:ug}, the prime triad can generate an octonion sized algebra, either $\OB$ or one of $\UB0$.
\end{definition}
\tab{ht}{aeg}{Non-Associativity to Associativity Examples}
\centering 
\begin{tabular}{|c|l|l|l|}  \hline   
{\bf Type} &{\bf\hfill Class~1\hfill} &{\bf\hfill Class~2\hfill} &{\bf\hfill Class~3\hfill} \\\hline
{\bf A} &$[\o1,\o{34},\o2] = -2\o{1234}$ &$[\o1,\o2,\o{34}]=0$ &$[\o2,\o1,\o{34}]=0$\\
{\bf B} &$[\o1,\o{345},\o{24}]=0$ &$[\o1,\o{24},\o{345}]=-2\o{1235}$ &$[\o{24},\o1,\o{345}]=0$\\
{\bf C} &$[\o1,\o{34},\o{24}]=0$ &$[\o1,\o{24},\o{34}]=0$ &$[\o{24},\o1,\o{34}]=-2\o{123}$\\
{\bf X} &$[\o1,\o4,\o2]=-2\o{124}$ &$[\o1,\o2,\o4]=2\o{124}$ &$[\o2,\o1,\o4]=-2\o{124}$ \\
\hline
\end{tabular}\tae

 The enumeration of 3-triad cycles in \tar{tss} can be divided by three and categorised as belonging to one of eight 3-triad cycles and the enumeration of non-cycles is separated into the four single triads rows shown in \tar{nas}. This table provides the 3-triad cycles structure for Cayley-Dickson algebras up to 10 generators and will be used to categorise zero divisors in the next section. As can be seen the structures up to $\UB4$ are distinct and after this the structure is the same so only $\UB{1}$, $\UB2$, $\UB3$ and $\UB4$ provide algebras with unique structure for ultracomplex numbers. The following theorem validates this table.

\tab{ht}{nas}{Non-Associative Structure}
\begin{tabular}{|c|c|c|c|c|c|c|c|c|c|}  \hline   
{\bf Structure}& $\OB$& $\UB1$& $\UB2$& $\UB3$& $\UB4$& $\UB5$& $\UB6$ &$\UB7$ \\\hline
{\bf AAA} &0 &28 &252 &1988 &14,868 &111,300 &849,492 &6,605,284 \\
{\bf BBA} &0 &0  &168 &1848 &15,120 &115,080 &872,928 &6,721,848 \\
{\bf ACC} &0 &0  &84  &1344 &12,852 &105,672 &834,876 &6,569,136 \\
{\bf XBB} &0 &0  &0   &672   &9072  &87,360  &752,976 &6,219,360 \\
{\bf BXC} &0 &0  &0   &168   &4536  &60,648  &621,432 &5,628,168 \\
{\bf CAB} &0 &0  &0   &0     &3024  &31,744  &577,584 &5,431,104 \\
{\bf CCX} &0 &28 &336 &2660 &18,648 &129,612 &931,392 &6,955,060 \\
{\bf XXX} &8 &64 &400 &2480 &16,368 &115,824 &863,600 &6,651,760 \\
{\bf A}   &0 &0  &84  &1008 &9828  &87,864  &747,180 &6,175,008 \\
{\bf B}   &0 &0  &0   &504   &7560  &78,456  &709,128 &6,022,296 \\
{\bf C}   &0 &28 &252 &1988 &14,868 &111,300 &849,492 &6,605,284 \\
{\bf X}   &4 &32 &284 &2080 &14,988 &111,000 &856,340 &6,588,272 \\ \hline
\end{tabular}\tae
\begin{theorem}[Cycles theorem]\label{thm:ct}
Of the 64 possible non-associativity 3-triad cycles only cases AAA, ACC, XBB, BBA, BXC, CAB, CCX and XXX exist. These will be called {\it silos}.
\end{theorem}
\begin{proof}
\tar{nas} shows the count of 3-Triad Cycles and Non-Cycles from \tar{tss} computed by code described in the last section. Starting with A-non-associativity for the first cycle then showing that the second and third cycles have Class~2 associativity shows that these are either A or C-non-associative. Here we use Class~2 associativity $[b,c,d]=0$ and Class~3 associativity $[c,b,d]=0$ and apply to the second cycle equivalent,
\begin{align}
  [b,bc,d] &= b(bc)d -b(bcd) = -bcbd +bcdb = cbbd +b(cd)b = -cd +cd = 0 \text{ and} \nonumber\\
  [c,bc,d] &= c(bc)d -c(bcd) = -bccd +bcdc = bd -c(bd)c = bd -bd = 0\es\label{eqn:tca}
\end{align}
Hence the result is either A or C and AA means third cycle A and AC means third cycle C needs to be shown. Using Class~3 associativity, $[c,b,d]=0$, then third cycle Class~2 associativity is
\[ [c,bc,d] = c(bc)d -c(bcd) = -bccd -cbdc = bd -c(bd)c = 0\ec \]
meaning either AAA or AAC. Now assume third cycle Class~3 associativity, then
\[ [bc,c,d] = bccd -bc(cd) = 0 \implies db = -bc(cd) \]
Substitute this into the second Class~3 associativity,
\[ [bc,b,d] = bcbd -bc(db) = -cbbd +bc(bc(cd)) = cd +cd(bc)(bc) = 0\ec \]
which means the third cycle for AA is also A. To the contrary, if the second cycle is not Class~3 associative then neither is the third so it must be ACC.

For the BBA and BXC silos we need to show that the second cycle has non-zero Class~2 associator and third cycle Class~2 associativity. Starting with the Class~2 associator for the second and third cycles and using Class~3 associativity, $[c,b,d]=0$, and Class~1 associativity, $[b,d,c]=0$, then
\begin{align*}
  [b,bc,d] &= b(bc)d -b(bcd) = -bcbd +bcdb \\
       &= cbbd -c(bd)b = dc +b(dc)b = 2dc \ne 0 \text{ and}\\
  [c,bc,d] &= c(bc)d -c(bcd) = -bccd +bcdc \\
       &= bd -c(bd)c = bd -bd = 0\es
\end{align*}
So the second cycle is either B or X-non-associative and the third A or C. If the second cycle is B then the third can be expanded using Class~3 and Class~1 associativity. Considering third cycle Class~1 associativity and using Class~3 associativity and assuming the third cycle Class~3 associativity from (\ref{eqn:tca}), then
\[ \begin{split}
    [c,d,bc] &= cd(bc) -c(d(bc)) = cd(bc) -bcdc = cd(bc) +c(bd)c \\
      &= cd(bd) +bc(cd) =cd(bd-bc) \ne 0\es
\end{split} \]
Hence BB has A as the third cycle and therefore BXC is the other 3-cycle.

For CAB and CCX silos, consider the second and third Class~2 associativity. Using Class~2 associativity, $[b,c,d]=0$, and Class~1 associativity, $[b,d,c]=0$, then
\begin{align*}
[b,bc,d] &= b(bc)d -b(bcd) = -bcbd +bcdb = cbbd +b(cd)b = -cd +cd = 0 \text{ and}\\
[c,bc,d] &= c(bc)d -c(bcd) = -bccd +bcdc = bd +b(cd)c = bd -bdcc = 2bd \ne 0\es
\end{align*}
Hence the second cycle is either A or C and the third either B or X. This is separated by showing that the second and third associator are either both associative or both non-zero. Substituting second Class~3 associativity, $[bc,b,d]=0$ as $cd = bc(db)$ into the third gives,
\[ [bc,c,d] = bccd -bc(cd) = -bd -bc(bc(db)) = -bd -db(bc)(bc) = 0\es \]
Hence C in the first cycle leads to CAB or CCX silos.

The silos that start with X non-associativity are proved using Moufang's non-associ\-ativity. From Theorem~\ref{thm:ma}, for X-non-associativity the Moufang types are 1 and 2,
\begin{align}
d(b(dc)) -dbdc = dcbd +bddc = 0 &\quad\Rightarrow\quad bc=dcbd\text{ and} \label{eqn:mt1}\\
b(d(cd)) -bdcd = cddb -bdcd = 0 &\quad\Rightarrow\quad bc=bdcd \label{eqn:mt2}\es
\end{align}
So the second cycle is either B or X. Class~3 and 1 for the second cycle using (\ref{eqn:mt2})  for Class~3 is
\begin{align}
[bc,b,d] &= bcbd-bc(bd)=cd-bc(bd)\ec \label{eqn:mt3}\\
[b,d,bc] &= bd(bc)-b(d(bc))=bd(bc)-bcdb \\
    &=bd(bc)-dcbddb=bd(bc)-dc = [bc,b,d]\es\nonumber
\end{align}
Class~2 and 3 for the third cycle with (\ref{eqn:mt1}) applied to Class~2 is
\begin{align}
[c,bc,d] &= c(bc)d -c(bcd) = bd +bcdc = bd +bdcddc = 2bd \ne 0\ec \nonumber \\
[bc,c,d] &= bccd -bc(cd) = -bd -bc(cd)\es \label{eqn:t3c2}
\end{align}
Hence, the third cycle is B or X, since Class~2 is non-zero, assuming the first cycle is X or B, from (\ref{eqn:mt1}). If Class~3 second cycle is zero, such that $cd = bc(bd)$ from (\ref{eqn:mt3}), then (\ref{eqn:t3c2}) becomes
\[ -bd -bc(cd) = -bd +bc(bd)(bc) = -bd -bd(bc)(bc) = 0\es \]
So again, the third cycle is B if the second is B, alternatively, both are X. Since the first cycle was assumed to be X we have 3-triads cycles XXX or XBB.
\end{proof}
It will be seen later that all A non-associativity triads generate zero divisors but evaluation finds that the three triads in the AAA triad-cycles do not satisfy the Moufang identities, yet all other A non-associativity types within the triad-cycles and most non-cycles triads satisfy Moufang's type 3 identity and not the Malcev identity, matching the association shown in Theorem~\ref{thm:ma}. Hence, the A, B, C and X non-associativity types seem to be more fundamental than the Moufang types and the triad-cycles are needed to distinguish the different types of loops. 

Half of the 3-triad cycle silos satisfy the identity for Malcev algebra, (\ref{eqn:mai}), which can now be recognised as second and third cycles from the Cycle theorem. Due to the cyclic rules the AAA and XXX silos satisfy (\ref{eqn:mai}) for all cycles so each 3-triads cycle provides a Malcev algebra of 3 triads. This is applicable to the eight XXX cycles of $\OB$ indicated in \tar{nas} and shown in (\ref{eqn:nox}). The four non-cycle X triads are $(\o2,\o3,\o{13})$, $(\o2,\o{13},\o{23})$, $(\o{12},\o3,\o{13})$ and $(\o{12},\o3,\o{23})$ and these also satisfy Malcev's identity since their second and third cycle parts overlap with other X triads. This completes the 28 Malcev triads in $\OB$ that provide a basis for a Malcev algebra. Since the scalars need to be avoided, the commutator is used or a product defined with zero norm, as provided in Elduque~\cite{Elduque}. There are four remaining silos with partial Malcev properties. The ACC and XBB satisfy (\ref{eqn:mai}) for the second and third cycles but not the first. The other two cases are BXC and CAB where the first cycle are Malcev algebra triads but the second and third are not. All these cases are derived from (\ref{eqn:mai}) which implies that only triad cycles that include A and X cases will satisfy this identity. This potentially identifies the basis of all Malcev algebras in the Cayley-Dickson algebras.

The X non-associative triad is the only completely symmetric non-associative type and this is the characteristic of octonions. This identifies the non-associative triads of $\OB$ in the XXX silo and the $\OB$ parts of $\UB1$ in \tar{nas}. Each triad generates 7 distinct blades giving $\binom73=28$ distinct triads and the members of 3-triad cycles are already included in this set so only the first cycle needs to be considered. Each octonion has $8\times3=24$ triads in the XXX silo and $4$ in the non-cycles. Algebra $\UB1$ has eight times these numbers and these can be identified with the following prime triads that generate the eight $\OB$ subalgebras, the first one expanded by (\ref{eqn:nox}),
\[ \begin{split}
  &(\o1, \o2, \o3), (\o1, \o2, \o4), (\o1, \o3, \o4), (\o1, \o{23}, \o4), \\
  &(\o2, \o3, \o4), (\o2, \o{13}, \o4), (\o{12}, \o3, \o4), (\o{12}, \o{13}, \o4)\es
\end{split} \]

The remaining AAA and CCX silos have 56 3-triad cycles giving seven copies of another algebra with prime triads
\begin{equation} \begin{split}
  &(\o1, \o2, \o{34}), (\o1, \o3, \o{24}), (\o1, \o{23}, \o{24}), (\o2, \o3, \o{14}), \\
  &(\o2, \o{13}, \o{14}), (\o{12}, \o3, \o{14}), (\o{12}, \o{13}, \o{14})\es 
  \end{split}\label{eqn:qop} \end{equation}
Each of these are isomorphic to an algebra documented in \cite{Wilmot} which is part of a series of six algebras that break the symmetry of octonions whereby subalgebras of exceptional Lie algebra $\GT$ provide automorphisms for these algebras. These have properties similar to the split octonions seen later, so were called the quasi-octonions, designated $\PB{k}$ where $k$ is the number of non-associative triads $k=4,8,10,12,14$ or $16$. The first case, $\PB4$, was discovered\cite{Cawagas1} and enumerated\cite{Cawagas2} by Cawagas. Only this algebra and two other quasi-octonions, $\PB{12}$ and $\PB{14}$, are subalgebras of this rank within ultracomplex numbers, the other three are not included. The number of occurrences of each subalgebra of the algebras shown in \tar{nas} are provided in \tar{usa}. Octonions have seven quaternion subalgebras, seen in (\ref{eqn:qa7}) for $(\o1,\o2,\o3)$, and the quasi-octonions have seven quaternion rings, split or not. Each quaternion-like algebra overlaps three times with other octonion-like loops. For example, $(\o1,\o2,\o{12})$ occurs in $(\o1,\o2,\o3)$, $(\o1,\o2,\o4)$ and $(\o1,\o2,\o{34})$. Hence the number of associative subalgebras in sedenions is $(7+8)\times7/3=35$ and from the table it can be seen that the number of overlapping quaternion-like algebras for each octonion-like subalgebra in $\UB{k}$ is $2^{k+1}-1$.

\tab{ht}{usa}{Ultracomplex Subalgebras}
\begin{tabular}{|c|c|c|c|c|c|c|c|c|c|}  \hline   
{\bf Structure}& $\OB$& $\UB1$& $\UB2$& $\UB3$& $\UB4$& $\UB5$& $\UB6$ &$\UB7$ \\\hline
$\HB$       &7 &35 &155 &651  &2667 &10,795 &43,435 &174,251 \\
$\OB$       &1 &8  &50  &310  &2046 &14,478  &107,950  &831,470 \\
$\PB4$     &0 &7  &63  &413  &2583 &16,905  &118,251  &873,901 \\
$\PB{12}$  &0 &0  &42  &504  &4158 &30,996  &229,194  &1,729,728 \\
$\PB{14}$  &0 &0  &0   &168  &3024 &34,776  &332,640  &2,912,616 \\\hline
Total $\PB{k}$ &0 &7 &105 &1085 &9765 &82,677  &680,085  &5,516,245 \\\hline
\end{tabular}\tae
\begin{definition}
   Define the number of  blades or pure basis elements in $\UB{m}$ as $N_{m+3}=2^{m+3}-1$. Also define $H_{m+3}$, $S_m$ and $O_m$ as the number of total $\HB$, $\PB{k}$ and $\OB$ copies in $\UB{m}$, respectively. Each $\PB{k}$ and $\OB$ consists of 28 unique non-associative triads, unique due to Theorem~\ref{thm:ug} for $\OB$. Hence the total number of triads in $\UB{m}$, which is the number of faces in the $(m+2)$-simplex, is 
\begin{equation}\label{eqn:OSm}
 \binom{N_{m+3}}3=H_{m+3} +28\times(O_m +S_m)\ec\qquad \text{$m\ge0$.}
\end{equation}
\end{definition}

\begin{theorem}[Ultracomplex algebra growth]\label{thm:u7}
 The number of $\PB{k}$ subalgebras in $\UB{m}$ is 
\begin{equation}\label{eqn:Om}
 S_{m+1} = 7\times(O_m +S_m)\ec\qquad\text{$m\ge0$.}
\end{equation}
\end{theorem}
\begin{proof}
Again, for convenience, allow $\UB0$ to include $\OB$. \tar{usa} shows these numbers for the Total $\PB{k}$ equal to $S_m$, for $m\in\N_0^7$. Including the scalar then Lemma~\ref{lem:1} shows that a product of the elements of a triad generates a fourth unique element. The set of all such triple products will already contain each fourth element so the total number of unique four-fold products is $\frac14\binom{2^{m+3}}3$. These are called quads but we need to remove products involving the scalar, which are produced by quaternions. Thus, 
\begin{equation}\label{eqn:Sm1}
 S_{m+1} = \frac14\binom{2^{m+3}}3 - H_{m+3}\es
\end{equation}
Note that this is $S_{m+1}$ not $S_m$ because power-associative algebras are generated at the next higher level. Multiplying each quaternion element in turn changes the quaternion triad into an X-non-associative type, which generates an octonion loop because it involves three generators. Multiplying each element of an octonion triad by a higher generator turns the X-non-associativity into types A, B and C, in turn. Further such doubling tends to keep to non-X types and, by Theorem \ref{thm:ug}, these will always be part of quasi-octonions because they involve more than three generators.

Substituting (\ref{eqn:OSm}) into (\ref{eqn:Om}) gives
\[ \begin{split}
    S_{m+1} &= \frac7{28}\Bigg(\binom{N_{m+3}}3-H_{m+3}\Bigg)\\
    &= \frac1{24}\Big((2^{m+3}-3)(2^{m+3}-1)(2^{m+3}-2)\Big) -\frac14H_{m+3}\\
    &= \frac14\Bigg(\binom{2^{m+3}}3 -3H_{m+3} -H_{m+3}\Bigg)\ec
\end{split}
\]
which matches (\ref{eqn:Sm1}), thus verifying (\ref{eqn:Om}). The factor of seven growth of the power-associative algebras is shown in the next section to be important for zero divisors. 
\end{proof}
\begin{corollary}
  The number of octonion and quasi-octonion algebras in $\UB{m}$, $T_m=O_m+S_m$, is 
  \[T_m=\frac1{21}N_{m+1}N_{m+2}N_{m+3}\ec\text{\qquad $m\ge0$.}\]
\end{corollary}
\begin{proof}
Byrapuram, et al.~\cite{Byrapuram}, derived (\ref{eqn:Sm1}) from a generalisation of the card game EvenQuads and expanded it to 
\begin{equation} \label{eqn:Sm}
  S_{m+1}=\frac1{24}\big(2^{m+3}-1)(2^{m+3}-2)(2^{m+3}-4\big)\es
\end{equation}
They noted that since powers of 2 can only have remainders 1, 2, or 4 modulo 7, then $S_{m+1}$ is divisible by 7, deriving the sequence A006096 (in the On-Line Encyclopedia of Integer Sequences, \url{https://oeis.org}, formerly M4982~\cite{Sloane}), so that (\ref{eqn:Om}) becomes
\begin{equation} \label{eqn:OS7}
  O_m+S_m = \frac1{21}(2^{m+3} -1)(2^{m+2} -1)(2^{m+1} -1)\ec
\end{equation}
which is easily confirmed. Substituting $N_m=2^m-1$ proves the theorem.
\end{proof}
\begin{theorem}[Ultracomplex subalgebras]\label{thm:us}
The number of Octonion subalgebras, $O_m$, and the sum of the number of quasi-octonion subalgebras, $S_m$, for Cayley-Dickson algebra $\UB{m}$ with $N_{m+3}=2^{m+3}-1$ pure basis elements are
\[ \begin{split}
    O_m &= \frac1{168} N_{m+1} N_{m+2} (N_{m+3} +49)\text{ and}\\
    S_m &= \frac1{24} N_{m+2} (N_{m+2} -1) (N_{m+2} -3)\es
\end{split}\]
\end{theorem}
\begin{proof}
The value for $S_m$ is gained by substituting $\N_{m+2}$ into (\ref{eqn:Sm}). Substituting this into (\ref{eqn:OS7}) gives the new sequence
\[ \begin{split}
  O_m &= \frac1{21}N_{m+3}N_{m+2}N_{m+1} 
        -\frac1{24}N_{m+2}(N_{m+2}-1)(N_{m+2}-3)\\
      &= N_{m+2}N_{m+1}\Big(\frac1{21}N_{m+3} 
          -\frac1{24}(N_{m+3}-7)\Big)\\
      &= \frac1{168}N_{m+2}N_{m+1}(N_{m+3}+49)\es
\end{split}\]
\end{proof}
Comparing the $O_m$ sequence using the $\OB$ row in \tar{usa} with the XXX row of \tar{nas} shows that octonion subalgebras of Cayley-Dickson algebras always have eight independent XXX triads.

The silo sources for each subalgebra are provided in \tar{sai} along with the non-associative type counts. This identification explains the structure shown in \tar{nas}. Each algebra has 8 3-triad cycles from the silos totalling 24 triads and the remaining 4 triads are always found in the non-cycles, so these are not identified but can be discovered by subtracting the identification column from the silos. It is interesting to note that all subalgebras have four X non-associativity types so a selection of a triad of type X does not guarantee an algebra isomorphic to $\OB$ will be generated.
\tab{ht}{sai}{Subalgebra Identification}
\begin{tabular}{|c|l|l|}  \hline 
{\bf Subalgebra} &{\bf Identification} &{\bf Silo counts} \\\hline
$\OB$     & 28X             & 8XXX \\
$\PB4$    & 12A, 12C, 4X    & 4AAA and 4CCX \\
$\PB{12}$ & 8A, 8B, 8C, 4X  & 4BBA, 2ACC, 2CCX or \\
          &                 & 2AAA, 2BBA, 2BXC, 2CCX \\
$\PB{14}$ & 7A, 10B, 7C, 4X & AAA, 3ACC, 4XBB or \\
          &                 & AAA, XBB, 3BXC, 3CAB or \\
          &                 & BBA, 2ACC, 2XBB, BXC, 2CAB or \\
          &                 & BBA, ACC, XBB, 2BXC, 3CAB \\\hline
\end{tabular}\tae

\begin{theorem}[Projection Theorem]\label{thm:pt}
Any triad from an ultracomplex algebra will generate an algebra isomorphic or anti-isomorphic to $\HB$, $\OB$, $\PB4$, $\PB{12}$ or $\PB{14}$. Furthermore, triads from $\UB1$ can project (or anti-project) to $\HB$, $\OB$ or $\PB4$, from $\UB2$ to a $\UB1$ subalgebra or $\PB{12}$ and all higher levels can project to $\UB2$ subalgebras or $\PB{14}$. This leads to the diagram of injections shown in Figure \ref{fig:coh}.
\begin{figure}[ht] \begin{center}
  \begin{tabular}{c}
  \includegraphics[width=10cm]{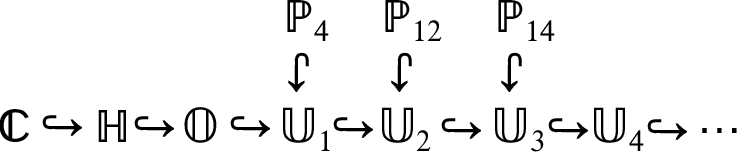}
  \end{tabular}
  \caption{\bf Injections of Cayley Dickson algebras} \label{fig:coh} 
\end{center} \toobig\ \end{figure} \toobig
\end{theorem}
\begin{proof}
The injections of the normed composition algebras, $\C$, $\HB$, $\OB$, are well defined by doubling, up to $\OB$, which has seven copies of $\HB$ and a triple non-associative part given by 28 triads. The former are rings and correspond to independent faces of the 6-simplex, as discussed in \cite{Wilmot}. The associative and non-associative parts can be isolated, with associative triads being rings by Theorem \ref{thm:at}, and the non-associative injection into the ultracomplex algebras analysed separately.
The number of distinct $\OB+\PB{k}$ subalgebra copies for each of the listed ultracomplex algebras matches the number of total non-associative triads divided by 28, as shown in \tar{zdall}, hence these subalgebra copies provide complete coverage of each ultracomplex algebra. A member of the XXX silo, $(\o1, \o2, \o3)$, generates $\OB$ and the first ordered member of CCX, $(\o1, \o{24}, \o{34})$, generates $\PB4$. \tar{nas} shows that these are the only silos available for $\UB1$ and \tar{sai} shows without B non-associativity only $\PB4$ can be included at this level. Then $\UB2$ contains these silos and BBA with the first ordered member, $(\o1, \o{24}, \o{345})$, generating $\PB{12}$. The 63 $\PB4$ algebras at this level use all the AAA and all but 84 CCX cycles, meaning the first identification in \tar{sai} provides the 42 $\PB{12}$. After this $\UB3$ contains silo XBB with first ordered member $(\o{14}, \o{25}, \o{36})$, which allows for $\PB{14}$ and the other identification of $\PB{12}$. Higher level ultracomplex algebras contain CAB, which can only generate $\PB{14}$, for example the prime, $(\o{14}, \o{256},\o{347})$, that first appears in $\UB4$.
\end{proof}
As shown later, each of the quasi-octonion algebras has 12 zero divisor pairs making them power-associative and non-alternative so this property flows through to the ultracomplex algebras. Since there are always seven copies of any of the quasi-octonion algebras then zero divisor pairs appear in multiples of 84.

\section{Zero Divisors}
\begin{definition} A {\it zero divisor} is a term whereby another term exists such that both are non-zero and the product is zero. Zero divisors can not be single blades and are not nilpotents of split algebras, so a {\it zero divisor pair} is defined, $(a+b)$ and $(c+d)$, such that $(a+b)(c+d) = 0$, for $a,b,c,d$ grade independent basis elements and $a^2=b^2$ and $c^2=d^2$. Linear combinations of these pairs are not considered.
\end{definition}
\begin{theorem}[Zero divisors theorem]\label{thm:zd}
The parallel terms and cross terms of the zero divisors equation separate into
\begin{equation}
  ac+bd = 0 \text{ and } ad+bc = 0\ec\label{eqn:zdt}
\end{equation}
with the second equation being an identity thus reducing the definition to just $ac=-bd$, as long as none of the terms form an associative algebra. Further, this is equivalent to Class~3 associativity and means the number of zero divisor pairs for $\UB{n}$, $Z_m$, $n=m+3>0$, and with $N_n=2^n-1$, is
\[ Z_m = \frac1{2}N_{m+2}(N_{m+2}-1)(N_{m+2}-3)\es \]
\end{theorem}
\begin{proof} Uniqueness of $a,b,c$ and $d$ ensures $ac$ does not have the same blade as either $ad$ or $bc$, from Lemma~\ref{lem:1}, so the two equations (\ref{eqn:zdt}) have distinct scaled blades and are separable. Since $b,c,d$ are distinct blades and $b,c,d$ are non-associative, then $a\approx bcd$, so $a$ is non-scalar or pure and all elements anticommute. Hence the equations can be re-arranged as
\begin{align}
  a = -(bd)c^{-1} = c(bd)/c^2 &\text{ and } a = -(bc)d^{-1} = cbd/d^2\text{ and }\label{eqn:zd1} \\
  d = acb^{-1} = b(ca)/b^2 &\text{ and } d = (bc)a^{-1} = bca/a^2\ec \label{eqn:zd2}
\end{align}
which implies $c^2=d^2$ and $a^2=b^2$. Substituting (\ref{eqn:zd1}) and (\ref{eqn:zd2}) shows Class~3 associativity
\[\begin{split}
 [b,c,a] &= bca -b(ca) = da^2 -db^2 = 0 \text{ and} \\
 [c,b,d] &= cbd -c(bd) = ad^2 -ac^2 = 0\es
\end{split}\]
For Type C non-associativity then $[b,d,c]=[b,c,d]=0$, which follows from Class~1 and~2 associativity. Assuming (\ref{eqn:zdt}) and multiplying the equations by $c$ and $d$, respectively, gives $ac^2+bdc = 0$ and $ad^2+bcd=0$ since $c^2=d^2$ then the $a$ terms cancel to give $bdc-bcd=0$. Applying the two associators gives $b(dc)-b(cd)=2b(dc)\ne0$, which is a contradiction. Hence Type~C non-associativity means $(a+b)(c+d)\ne0$. Type~X also does not generate zero divisors leaving Class~3 associativity. Note that the reversed product $(c+d)(a+b)$ generates $-ac-bd=0$ so provides the same zero divisor pair.
Finally, substituting for $a$, $ad = cbdd/d^2 = -bc$ so that $ad+bc=0$ is an identity. Moreno~\cite{Moreno} first showed that $d=-acb$ and  Zhilina~\cite{Zhilina} extended this work to show that the zero divisor equation reduces to one of (\ref{eqn:zdt}). These derivations implicitly include the ordering $a<b$ and $c<d$ provided by the usual Cayley-Dickson product $(a,b)(c,d)$ when taking the $\o{n}$ generators into account but neither author included that Type~C excludes zero divisors.

Hence zero divisor pairs are separated into two types for A and B non-associative triads that are both Class~3 associative. Octonions can not have zero divisors because the non-associative triads are X non-associative. Swapping the last two elements of the associator in \tar{tat} swaps Class~1 associativity with Class~2 associativity and has no correspondence with Class~3 associativity other than changing the sign. Hence Type~A and~B zero divisors may overlap because swapping the terms within brackets produces the same zero divisor pair. \tar{sai} shows 12 A non-associative types for $\PB4$, which corresponds to 12 zero divisor pairs. The other quasi-octonion algebras have mixtures of A and B types adding up to more than 12 but since all ordered triples are included in the list of 28 triads then some A and B types must overlap as zero divisors to give 12 distinct zero divisor pairs for each quasi-octonion, which was found in~\cite{Wilmot} by enumerating all triads in each product table. Theorem~\ref{thm:us} defines the number of quasi-octonion algebras at each level as $S_m=\frac1{24} N_{m+2} (N_{m+2} -1) (N_{m+2} -3)$.
Since each $\PB{k}$ algebra has 12 zero divisor pairs~\cite{Wilmot} then $Z_m=12S_m$.
\end{proof}
Moreno~\cite{Moreno} defines a triple in ${\rm A}_n$ as special if $[a,y,z]\ne0$ for $a, z$ distinct, pure basis elements. For zero divisors of the form $(a,z)(c,y)=(0,0)$ then the associator translates to $[a,d,b]$ where $b=z\o{n}$ and $d=y\o{n}$. Multiplying both $b$ and $d$ by $\o{n}$ then by Lemma~\ref{lem:4} returns the Moreno associator so that $[a,d,b]\ne0$. From \tar{tat} this is Class~3 non-associativity. Moreno's Theorem~2.13 proves this triple is octonion, which matches the Type~X non-associativity of the associator. Theorem~2.9 then assumes $c=-(ay)z$ to find this satisfies the zero divisor equation and the zero divisor pairs listed below in \tar{zds} satisfy the translated form $c=-adb$. Now, Class~3 non-associativity also satisfies Type~C non-associativity but these do not define zero divisors. For example, the $\SB$ equation $adb=\o{123}\o{34}\o1=-\o{24}=-c$ and yet $(\o{123}+\o1)(\o{24}+\o{34})=-2\o{134}$. Zhilina~\cite{Zhilina} derives this same condition and for $(a,z)(c,y)=0$, assuming pure basis elements then Lemma 3.14 has $(c,-(zc)a)=(-(ay)z,y)$. For norm $+1$ and conjugate of $y$ being $-y$, then $c=-ayz$ and $-zca=y$ or $-zc-ay=0$ and $zc+ay=0$, which is one of the zero divisor statements in Theorem~\ref{thm:zd}. The remark 3.15 for $\chi=1$ has $[z,c,a]=0$ and $[a,y,z]=0$ which are both Class~3 associativity implying Type~A or Type~B non-associativity. Hence this is a more complete exposition of zero divisors. Biss~\cite{Biss} and Cohen~\cite{Cohen} derive the same equations but do not mention the two associativities needed to restrict the zero divisor conditions.

Since the quasi-octonions, $\UB0$, contain zero divisors that can't be expressed using just loop multiplication, the designation non-associative ring is more appropriate, as suggested by Cohen~\cite{Cohen}. The term subalgebra of the sedenion algebra is also appropriate and the terminology ultracomplex distinguishes $\UB0$ from quasi-octonions that can't be represented as products of Cayley-Dickson elements. Only $\OB$ and $\PB4$ form non-associative loops within sedenions but \tar{sai} shows there are many representations of other $\UB0$ algebras within higher level Cayley-Dickson algebras. The cycles of Theorem~\ref{thm:ct} provide a new way to analyse these ultracomplex numbers and we now consider the structure within each non-associative ring and for zero divisors.
\begin{definition}[Modes] Calling (\ref{eqn:zdt}) the {\it prime zero divisor pair} then $(-d+b)(c+a)$, is called the {\it dual of the prime}. The terminology dual is appropriate for zero divisors. Also $(a^\prime+b)(c+\abs(db))$, with $a^\prime=bc \abs(db)$ and $\abs(db)$ being the absolute value of $db$, is called the {\it extended zero divisor pair}, and  $(-\abs(db)+b)(c+a^\prime)$, is called {\it both extended and dual}. These four cases are called modes and whether they are zero divisors determines another level of structure for the ultracomplex algebras. 
\end{definition}
The definition of zero divisor pairs required that $a^2=b^2$ and $c^2=d^2$. For the extended equation to hold these become $a^{\prime2}=b^2$ and $c^2 = (db)^2$, which implies $b^2=c^2=-1$. Hence for all subsequent work we assume $a^2 = b^2 = c^2 = d^2 = \pm1$ and take $b$, $c$ and $d$ to be positive blades.
\begin{theorem}[Modes theorem]\label{thm:mt}
Triad cycle mode transformations behave as described in the following table.
\begin{center}
\begin{tabular}{|c|c|c|c|}  \hline 
{\bf Primary} & {\bf Dual} & {\bf Extended} & {\bf Both}  \\ \hline
 AAA & AAA & AAA & AAA \\
 BBA & CCX & BBA & CCX \\
 ACC & ACC & XBB & XBB \\
 XBB & XBB & ACC & ACC \\
 BXC & CAB & CAB & BXC \\
 CAB & BXC & BXC & CAB \\
 CCX & BBA & CCX & BBA \\
 XXX & XXX & XXX & XXX \\\hline
\end{tabular}
\end{center}
\end{theorem}
\begin{proof} 
The dual mode replaces $d$ with $a$ in the triad $(b,c,d)$ so substituting this into \tar{tat} we find,
\begin{align*}
\text{Class~1 associativity }[b,d,c] &\rightarrow [b,a,c] \text{ which is Class~1 associativity,} \\
\text{Class~2 associativity }[b,c,d] &\rightarrow [b,c,a] \text{ which is Class~3 associativity (reversed),} \\
\text{Class~3 associativity }[c,b,d] &\rightarrow [c,b,a] \text{ which is Class~2 associativity (reversed).}
\end{align*}
Together these mean that under the dual transformation for the first cycle B $\leftrightarrow$ C and A and B are independent. The Modes Theorem table shows that A $\rightarrow$ A and  X $\rightarrow$ X for the first primary triad, and either this case or A $\leftrightarrow$ X for the other cycles. This can occur because the dual transformations above allow Class~2 and~3 associativities to both change for the A and X types. These cases need to be analysed in detail.

For traceability, the terminology ${\rm P}p_q$ is introduced to identify a prime Class~$p$ assoc\-iativity equation for the $q^{\rm th}$ cycle, where $p,q\in\N_1^3$ and ${\rm D}p_q$ for the equivalent classes and cycles for the dual mode. The following class, cycles and mode associativity equations need to be evaluated 
\begin{align}
[b,d,c] &= bdc -b(dc) =bdc +dcb\ec\tag{${\rm P}1_1$}\label{eqn:p11} \\
[b,c,d] &= bcd -b(cd) =bcd +cdb\ec\tag{${\rm P}2_1$}\label{eqn:p21} \\
[c,b,d] &= cbd -c(bd) =cbd +bdc =-bcd -dbc\ec\tag{${\rm P}3_1$}\label{eqn:p31} \\
[b,d,bc] &= bd(bc) -b(d(bc)) =bd(bc) -bcdb\ec\tag{${\rm P}1_2$}\label{eqn:p12} \\
[bc,b,d] &= bcbd -bc(bd) =cd -bc(bd)\ec\tag{${\rm P}3_2$}\label{eqn:p32} \\
[c,a,bc] &= ca(bc) -c(a(bc)) =-ac(bc) +a(bc)c\ec\tag{${\rm D}1_3$} \label{eqn:d13} \\
[b,bc,a] &= b(bc)a -b(bca) =ac -a(bc)b\ec\tag{${\rm D}2_2$}\label{eqn:d22} \\
[c,bc,a] &= c(bc)a -c(bca) =-ab -a(bc)c\es\tag{${\rm D}2_3$}\label{eqn:d23}
\end{align}
Since A$\rightarrow$A or A$\rightarrow$X then we only need to find Class~2 or Class~3 associativity. Firstly, note that ${\rm P}3_1=0$ implies primary first cycle A or B and since, by definition, $a=bcd$ then $a=dbc$. Substituting this and $a=bcd=-d(bc)$ into ${\rm D}1_3$ gives
\[{\rm D}1_3=-ac(bc)+a(bc)c=dbcc(bc)+d(bc)(bc)c=-bc(bd)+cd\es\]
Hence ${\rm D}1_3=\neg({\rm P}3_2\veebar{\rm P}3_1)$ so if the primary third cycle is A or X then the dual third cycle is A or X and via versa X or B imply the dual ends with C or B. Now consider ${\rm D}2_2$ and ${\rm D}2_3$ assuming ${\rm P}3_1=0$ in the first, so that $a=bdc$, and ${\rm P}2_1=0$, with $a=-cdb$, in the second and $a=-d(bc)$ in both then
\[ \begin{split}
  {\rm D}2_2 &= ac-a(bc)b=bdcc+d(bc)(bc)b =-bd-db=0\text{ and}\\
  {\rm D}2_3 &= -ab-a(bc)c =cdbb +d(bc)(bc)c =-cd -dc =0\es
\end{split} \]
Including the assumptions we have ${\rm D}2_2={\rm P}3_1$ and ${\rm D}2_3={\rm P}2_1$. Given the uniqueness of the silos this is enough information to define the dual mode, as shown in \tar{dmt}.
\tab{ht}{dmt}{Dual Mode Transformations}
\begin{tabular}{|c|c|c|c|c||c|c|c|c|c|}  \hline   
 \multicolumn{5}{|c||}{\bf Prime Mode} & \multicolumn{4}{|c|}{\bf Dual Mode} \\\hline
{\bf Silo} &${\rm P}2_1$ &${\rm P}3_1$ &${\rm P}3_2$ &$\neg({\rm P}3_1\veebar{\rm P}3_2)$
  &${\rm D}1_3$ &${\rm D}2_2$ &${\rm D}2_3$ &{\bf Silo}\\\hline
{\bf AAA} &0   &0   &0   &\nz  &\nz &0   &0    &{\bf AAA} \\
{\bf BBA} &\nz &0   &0   &\nz  &\nz &0   &\nz  &{\bf CCX} \\
{\bf ACC} &0   &0   &\nz &0    &0   &0   &0    &{\bf ACC} \\
{\bf XBB} &\nz &\nz &0   &0    &0   &\nz &\nz  &{\bf XBB} \\
{\bf BXC} &\nz &0   &\nz &0    &0   &0   &\nz  &{\bf CAB} \\
{\bf CAB} &0   &\nz &0   &0    &0   &\nz &0    &{\bf BXC} \\
{\bf CCX} &0   &\nz &\nz &\nz  &\nz &\nz &0    &{\bf BBA} \\
{\bf XXX} &\nz &\nz &\nz &\nz  &\nz &\nz &\nz  &{\bf XXX} \\ \hline
\end{tabular}\tae

Using identifiers ${\rm E}p_q$ for the $p$ cycle and $q$ class extended mode associativity equations then we require Classes~1 and 3 for the first and second cycles and Classes~2 and 3 for the third cycle,
\begin{align}
[b,\abs(db),c]&= \pm(b(db)c-b(dbc))=\pm(dc-bdcb)\ec\tag{${\rm E}1_1$}\label{eqn:e11} \\
[c,b,\abs(db)]&=\pm(cb(db)-c(b(db)))=\pm(bc(bd)-cd)\ec\tag{${\rm E}3_1$}\label{eqn:e31} \\
[b,\abs(db),bc]&= \pm(b(db)(bc)-b(db(bc)))=\pm(-bcd-bd(bc)b)\ec\tag{${\rm  E}1_2$}\label{eqn:e12} \\
[bc,b,\abs(db)]&=\pm(bcb(db)-bc(b(db)))=\pm(-dbc-bcd)\ec\tag{${\rm E}3_2$}\label{eqn:e32} \\
[c,bc,\abs(db)]&=\pm(c(bc)(db)-c(bc(db)))=\pm(d-bc(bd)c)\ec\tag{${\rm E}2_3$}\label{eqn:e23}  \\
[bc,c,\abs(db)]&=\pm(bcc(db)-bc(c(db)))=\pm(-d+bc(dbc))\es\tag{${\rm E}3_3$}\label{eqn:e33}
\end{align}

\tab{ht}{emt}{Extended Mode Transformations}
\begin{tabular}{|c|c|c|c|c||c|c|c|c|c|c|c|c|}  \hline   
 \multicolumn{5}{|c||}{\bf Prime Mode} & \multicolumn{7}{|c|}{\bf Extended Mode} \\\hline
{\bf Silo} &${\rm P}1_1$ &${\rm P}3_1$ &${\rm P}1_2$ &${\rm P}3_2$
  &${\rm E}1_1$ &${\rm E}3_1$ &${\rm E}1_2$ &${\rm E}3_2$ &${\rm E}2_3$ &${\rm E}3_3$ &{\bf Silo}\\\hline
{\bf AAA} &\nz &0   &\nz &0    &\nz &0   &\nz &0   &0   &0   &{\bf AAA} \\
{\bf BBA} &0   &0   &0   &0    &0   &0   &0   &0   &0   &0   &{\bf BBA} \\
{\bf ACC} &\nz &0   &0   &\nz  &\nz &\nz &0   &0   &\nz &0   &{\bf XBB} \\
{\bf XBB} &\nz &\nz &0   &0    &\nz &0   &0   &\nz &0   &\nz &{\bf ACC} \\
{\bf BXC} &0   &0   &\nz &\nz  &0   &\nz &\nz &0   &\nz &0   &{\bf CAB} \\
{\bf CAB} &0   &\nz &\nz &0    &0   &0   &\nz &\nz &0   &\nz &{\bf BXC} \\
{\bf CCX} &0   &\nz &0   &\nz  &0   &\nz &0   &\nz &\nz &\nz &{\bf CCX} \\
{\bf XXX} &\nz &\nz &\nz &\nz  &\nz &\nz &\nz &\nz &\nz &\nz &{\bf XXX} \\ \hline
\end{tabular}\tae
The following equivalences can be identified. ${\rm E}1_1\equiv{\rm P}1_1 b$, ${\rm E}3_1\equiv{\rm P}3_2$, ${\rm E}1_2\equiv{\rm P}1_2 b$, ${\rm E}3_2\equiv{\rm P}3_1$, ${\rm E}2_3\equiv{\rm P}3_2 c$ and ${\rm E}3_3\equiv{\rm P}3_1$. 
The later one is because if (\ref{eqn:p31}) is zero then $bc(dbc)=(bcd)(bc)=-d(bc)(bc)=d$ so that (\ref{eqn:e33}) is zero or both are non-zero.
This is used in \tar{emt} to transform all silos to the extended mode equivalents where \nz\  means non-zero. This proves the extended results and the both mode follows by transforming the dual mode as if it was an extended mode to generate the last column of the theorem. It is observed that the non-cycle triads transform according to the cycle $(b,c)$ with $d < bc$. Hence they behave like the first element of each primary triad with the second and third cycles occurring elsewhere.
\end{proof}
\begin{corollary}
  There are 12 zero divisors for each power-associative ring.
\end{corollary}
\begin{proof}
  The modes define the structure of the power-associative rings. Each mode of the triad cycles makes up twelve of the loop's triads so that two triad cycles and four non-cycle triads are required to make up the 28 non-associative triads of any $\UB0$ or $\OB$. The modes for any type B primary is seen in Theorem~\ref{thm:mt} to be duplicated half of the time and otherwise become type C over the other modes while the type A primaries are duplicated or change to two type B modes. Hence counting type A non-associativity as a single zero divisor and type B as half a zero divisor then the identification column from \tar{sai} provides the result; the 12A types of $\PB4$, the 8A plus 8B of $\PB{12}$ and the 7A plus 10B of $\PB{14}$ generate twelve zero divisors each. Note that the primary triads in \tar{nas} do not show all modes from Theorem~\ref{thm:mt} because some to the $d$ values are smaller than the $a$ or $b$ values, hence will appear elsewhere in this ordered table. The Silo Counts of \tar{sai} do show all representations provided by all the modes. It should also be noted that the silos of \tar{nas} provide a method to count zero divisors. Taking each A as one and each B as a half in the first column and multiplying the values in the algebra columns and summing gives the total zero divisors. For example, the CAB contribution for $\UB4$ is $1.5\times3024=4536$ zero divisors. 
\end{proof}

All zero divisor pairs in \tar{nas} are derived only from the silos that contain A and B non-associative types. The mode transformations for the AAA silo provides three other AAA 3-cycle triads. Hence the AAA silo can be reduced by factors of 3 cycles and 4 modes in all ultracomplex algebras. The BBA silo can also be reduced by 3 cycles but the modes are more complicated. Only the A dual mode and B extended mode generate zero divisors so modes provide a reduction by two only. It is observed that for the ultracomplex algebras shown in \tar{nas}, that for $\UB1\equiv\SB$ only the AAA silo can provide distinct zero divisor pairs, numbering $28\times3=84$.

Cawagas~\cite{Cawagas1} labels the pure sedenions as $\e{i}$, for $i\in\N_1^{15}$. The graded single, pure elements, called the basis of $\SB$, have a natural order generated by the Cayley-Dickson process, which match these in order as
\begin{equation}\begin{split} (\e1,\e2,\dots\e{15}) = 
  (&\o1, \o2, \o{12}, \o3, \o{13}, \o{23}, \o{123}, \o4, \\
   &\o{14}, \o{24}, \o{124}, \o{34}, \o{134}, \o{234}, \o{1234})\es
   \label{eqn:basis}\end{split}\end{equation}
The first three elements are recognised as the pure quaternions and the first seven match the octonions~\cite{Baez}. Cawagas derives the $84$ zero divisor pairs shown in \cite{Cawagas1} and these are provided in graded form in \tar{zds} as $b$ and $c$ columns followed by prime, dual, extended and its dual triads as $d$, $a$, $db$ and $ab$ columns, respectively. 

Sedenions are power-associative which means the associator with repeated terms is always zero, $[a,a,a]=0$. This is obvious for an algebra where the single terms square to a scalar and all terms have inverses. Octonions are alternate associative which means swapping terms of the associator changes the sign or $[a,a,b] = [a,b,a] = 0$ for all $a,b\in\OB$. But sedenions are not alternate associative and demonstrate unordered or asymmetric associativity described above. For example, the A type from \tar{aeg} is
\begin{align*}
\left[\o1,\o2,\o{34}\right] &= -\o{12}\o{34} -\o1\o{234} = 0 \text{ and }\\
\left[\o1,\o{34},\o2\right] &= \o{134}\o2 +\o1\o{234} = -2\o{1234}\es
\end{align*}
\tab{ht}{zds}{$\UB1$ Zero Divisors}
\begin{tabular}{|c|c|c||c|c|c|c|}  \hline   
  &{\bf\em b} &{\bf\em c} &{\bf\em d} &{\bf\em a}  &{\bf\em db}&{\bf\em ab} \\ \hline
1  &$\o{1}$   &$\o{2}$    &$\o{34}$   &$-\o{1234}$ &$\o{134}$  &$\o{234}$   \\
2  &$\o{1}$   &$\o{12}$   &$\o{34}$   &$\o{234}$   &$\o{134}$  &$\o{1234}$  \\
3  &$\o{1}$   &$\o{3}$    &$\o{24}$   &$\o{1234}$  &$\o{124}$  &$-\o{234}$  \\
4  &$\o{1}$   &$\o{13}$   &$\o{24}$   &$-\o{234}$  &$\o{124}$  &$-\o{1234}$ \\
5  &$\o{1}$   &$\o{23}$   &$\o{24}$   &$\o{134}$   &$\o{124}$  &$\o{34}$    \\
6  &$\o{1}$   &$\o{123}$  &$\o{24}$   &$-\o{34}$   &$\o{124}$  &$\o{134}$   \\
7  &$\o{2}$   &$\o{12}$   &$\o{34}$   &$-\o{134}$  &$\o{234}$  &$\o{1234}$  \\
8  &$\o{2}$   &$\o{3}$    &$\o{14}$   &$-\o{1234}$ &$-\o{124}$ &$-\o{134}$  \\
9  &$\o{2}$   &$\o{13}$   &$\o{14}$   &$\o{234}$   &$-\o{124}$ &$\o{34}$    \\
10 &$\o{2}$   &$\o{23}$   &$\o{14}$   &$-\o{134}$  &$-\o{124}$ &$\o{1234}$  \\
11 &$\o{2}$   &$\o{123}$  &$\o{14}$   &$\o{34}$    &$-\o{124}$ &$-\o{234}$  \\
12 &$\o{12}$  &$\o{3}$    &$\o{14}$   &$\o{234}$   &$\o{24}$   &$-\o{134}$  \\
13 &$\o{12}$  &$\o{13}$   &$\o{14}$   &$\o{1234}$  &$\o{24}$   &$\o{34}$    \\
14 &$\o{12}$  &$\o{23}$   &$\o{14}$   &$-\o{34}$   &$\o{24}$   &$\o{1234}$  \\
15 &$\o{12}$  &$\o{123}$  &$\o{14}$   &$-\o{134}$  &$\o{24}$   &$-\o{234}$  \\
16 &$\o{3}$   &$\o{13}$   &$\o{24}$   &$-\o{124}$  &$-\o{234}$ &$-\o{1234}$ \\
17 &$\o{3}$   &$\o{23}$   &$\o{14}$   &$\o{124}$   &$-\o{134}$ &$\o{1234}$  \\
18 &$\o{3}$   &$\o{123}$  &$\o{14}$   &$-\o{24}$   &$-\o{134}$ &$-\o{234}$  \\
19 &$\o{13}$  &$\o{23}$   &$\o{14}$   &$\o{24}$    &$\o{34}$   &$\o{1234}$  \\
20 &$\o{13}$  &$\o{123}$  &$\o{14}$   &$\o{124}$   &$\o{34}$   &$-\o{234}$  \\
21 &$\o{23}$  &$\o{123}$  &$\o{24}$   &$\o{124}$   &$\o{34}$   &$\o{134}$   \\\hline
\end{tabular}\tae
This example produces the zero divisor pair $(\o1-\o{1234})(\o2+\o{34})$ that provides an example of non-alternate associativity, $[a,a,b]\ne0$ even though $[a,b,a]=0$, 
\begin{align*}
\left[(\o1-\o{1234}),(\o1-\o{1234}),(\o2+\o{34})\right] &= -2\o2 -2\o{34} \text{ and }  \\
\left[(\o1-\o{1234}),(\o2+\o{34}),(\o1-\o{1234})\right] &= 0\es
\end{align*}
This is because each term is invertible so that $(\o1-\o{1234})^2 = -2$ but multiplying the zero divisor terms first gives zero. For $\UB1$, we found above that zero divisors need $d=\o{i4}$ or $d=\o{ij4}$ or $d=\o{1234}$, for $i,j\in\N_1^3$, but $b$ and $c$ do not contain $\o4$ in order to satisfy $[c,b,d]=0$. This means $a$ also contains the $\o4$ generator, satisfying the symmetry required for the dual and other modes. 

\tab{ht}{zd7}{AAA Primary Zero Divisors in $\UB1$}
\begin{tabular}{|c |c|}  \hline 
 {\bf\em bc}    &{\bf\em d} \\\hline
 $\o1\o2$       &$\o{34}$  \\
 $\o1\o3$       &$\o{24}$  \\
 $\o1\o{23}$    &$\o{24}$  \\
 $\o2\o3$       &$\o{14}$  \\
 $\o2\o{13}$    &$\o{14}$  \\
 $\o{12}\o3$    &$\o{14}$  \\
 $\o{12}\o{13}$ &$\o{14}$  \\\hline
\end{tabular}\tae
Sedenions contain $455$ ordered triples split into 35 associative and, from \tar{nas}, $84$ A non-associative, $84$ C non-associative and $252$ X or completely non-associative. Cawagas's 84 zero divisor pairs are summarised in AAA 3-triad cycle form in \tar{zd7} which shows the primes only. Applying the Mode and Cycles Theorems gives twelve zero divisor pairs, applicable for the AAA silo, over the seven 3-triad cycles generating all $84$ zero divisor pairs in \tar{zds}. Multiplying all combinations of the three elements in the rows of \tar{zd7} generate seven pure elements and these correspond to the seven $\PB4$ subalgebras of sedenions shown in (\ref{eqn:qop}) and modes and cycles correspond to other triads of each $\PB4$ subalgebra. Thus, the cycles and modes are a new way to analyse ultracomplex numbers.
\tab{ht}{zd12o}{AAA Primary Zero Divisors in $\UB2$ based on $\OB$}
\begin{tabular}{|c c c| c |c|c|}  \hline 
 &{\bf\em bc cycle}& &{\bf\em d}  &{\bf\em d$_1$} &{\bf\em d$_2$}\\\hline
 $\o1\o2$      &$\o1\o{12}$   &$\o2\o{12}$     &$\o{34}$  &$\o{35}$ &$\o{45}$\\
 $\o1\o3$      &$\o1\o{13}$   &$\o3\o{13}$     &$\o{24}$  &$\o{25}$ &$\o{45}$\\
 $\o1\o{23}$   &$\o1\o{123}$  &$\o{23}\o{123}$ &$\o{24}$  &$\o{25}$ &$\o{45}$\\
 $\o2\o3$      &$\o2\o{23}$   &$\o3\o{23}$     &$\o{14}$  &$\o{15}$ &$\o{45}$\\
 $\o2\o{13}$   &$\o2\o{123}$  &$\o{13}\o{123}$ &$\o{14}$  &$\o{15}$ &$\o{45}$\\
 $\o{12}\o3$   &$\o{12}\o{123}$ &$\o3\o{123}$   &$\o{14}$  &$\o{15}$ &$\o{45}$\\
 $\o{12}\o{13}$&$\o{12}\o{23}$&$\o{13}\o{23}$  &$\o{14}$  &$\o{15}$ &$\o{45}$\\
 \hline
\end{tabular}\tae

The Cayley-Dickson algebra with five generators, $\UB2$, has $4495$ ordered triples split into $155$ associative, $1092$ A, $336$ B, $1092$ C and $1820$ X non-associative triads. This matches the assignment from \tar{usa}, $\UB2=50\times\OB+63\times\PB4+42\times\PB{12}$.  Now $\UB{2}$ has $1260$ zero divisor pairs because unlike $\UB1$ with only AAA 3-triad cycles, \tar{nas} shows that BBA and ACC cases are involved. The AAA and BBA silos provide all the zero divisors so that the $84$ A types from the first cycle of the ACC silo and the $84$ non-cycle A types must overlap with the first two silos. So only the first two silos are analysed. What is not documented is that some of the BBA and AAA cases overlap, in that the same $b$, $c$ 3-cycle can be included in both silos for different $d$ elements. All these cases will be shown in different tables starting with ones based on octonion 3-cycles.
\tab{ht}{zd12}{AAA Primary Zero Divisors in $\UB2$ independent of $\OB$}
\begin{tabular}{|c c c|c|c|c|}  \hline 
 &{\bf\em bc cycle}& &{\bf\em d$_1$} &{\bf\em d$_2$} &{\bf\em d$_3$}\\\hline
  $\o{1}\o{4}$   &$\o{1}\o{14}$     &$\o{4}\o{14}$    &$\o{25}$  &$\o{35}$ &$\o{235}$\\
  $\o{2}\o{4}$   &$\o{2}\o{24}$     &$\o{4}\o{24}$    &$\o{15}$  &$\o{35}$ &$\o{135}$\\
  $\o{12}\o{4}$  &$\o{12}\o{124}$   &$\o{4}\o{124}$   &$\o{15}$  &$\o{35}$ &$\o{135}$\\
  $\o{3}\o{4}$   &$\o{3}\o{34}$     &$\o{4}\o{34}$    &$\o{15}$  &$\o{25}$ &$\o{125}$\\
  $\o{13}\o{4}$  &$\o{13}\o{134}$   &$\o{4}\o{134}$   &$\o{15}$  &$\o{25}$ &$\o{125}$\\
  $\o{23}\o{4}$  &$\o{23}\o{234}$   &$\o{4}\o{234}$   &$\o{15}$  &$\o{25}$ &$\o{125}$\\
  $\o{123}\o{4}$ &$\o{123}\o{1234}$ &$\o{4}\o{1234}$  &$\o{15}$  &$\o{25}$ &$\o{125}$\\
 \hline
\end{tabular}\tae
\tar{zd12o} shows solutions derived from $\UB1$ which used the 4\textsuperscript{th} generator in the $d$ term. It uses the Cycle and Mode Theorems to simplify the 84 zero divisor pairs of $\UB1$ to just 7 primaries and since they are AAA non-associative then all four modes can be derived from column $d$. The $bc$ 3-cycles cover all octonion pairs and as mentioned above only the $d$ term contains $\o4$. This is extended by $\UB2$ to include the 5\textsuperscript{th} generator then both together, $\o{45}$, as shown in the $d_1$ and $d_2$ columns of the table. All of these triads are AAA non-associative with 3-cycle solutions thus giving three lots of $84$ or $252$ zero divisor pairs. The remaining AAA cases not based on octonion 3-cycles are provided in \tar{zd12}.
\tab{ht}{zd8}{Overlapping AAA and BBA Primary Zero Divisors in $\UB2$}
\begin{tabular}{|ccc|c|c|c|}  \hline 
 &{\bf\em bc cycle}& &{\bf\em d} &{\bf\em d$_1$} &{\bf\em d$_2$}\\\hline
$\o{1}\o{24}$ &$\o{1}\o{124}$ &$\o{24}\o{124}$     &$\o{25}$  &$\o{345}$  &$\o{2345}$ \\
$\o{1}\o{34}$ &$\o{1}\o{134}$ &$\o{34}\o{134}$     &$\o{35}$  &$\o{245}$  &$\o{2345}$ \\
$\o{1}\o{234}$ &$\o{1}\o{1234}$ &$\o{234}\o{1234}$ &$\o{235}$ &$\o{245}$  &$\o{345}$ \\
$\o{2}\o{14}$ &$\o{2}\o{124}$ &$\o{14}\o{124}$     &$\o{15}$  &$\o{345}$  &$\o{1345}$\\
$\o{2}\o{34}$ &$\o{2}\o{234}$ &$\o{34}\o{234}$     &$\o{35}$  &$\o{145}$  &$\o{1345}$ \\
$\o{2}\o{134}$ &$\o{2}\o{1234}$ &$\o{134}\o{1234}$ &$\o{135}$ &$\o{145}$  &$\o{345}$\\
$\o{12}\o{14}$ &$\o{12}\o{24}$ &$\o{14}\o{24}$     &$\o{15}$  &$\o{345}$  &$\o{1345}$ \\
$\o{12}\o{34}$ &$\o{12}\o{1234}$ &$\o{34}\o{1234}$ &$\o{35}$  &$\o{145}$ &$\o{1345}$\\
$\o{12}\o{134}$ &$\o{12}\o{234}$ &$\o{134}\o{234}$ &$\o{135}$ &$\o{1245}$ &$\o{345}$ \\
$\o{3}\o{14}$ &$\o{3}\o{134}$ &$\o{14}\o{134}$     &$\o{15}$  &$\o{245}$ &$\o{2345}$\\
$\o{3}\o{24}$ &$\o{3}\o{234}$ &$\o{24}\o{234}$     &$\o{25}$  &$\o{145}$ &$\o{1245}$ \\
$\o{3}\o{124}$ &$\o{3}\o{1234}$ &$\o{124}\o{1234}$ &$\o{125}$ &$\o{145}$  &$\o{245}$\\
$\o{13}\o{14}$ &$\o{13}\o{34}$ &$\o{14}\o{34}$     &$\o{15}$  &$\o{245}$ &$\o{1245}$ \\
$\o{13}\o{24}$ &$\o{13}\o{1234}$ &$\o{24}\o{1234}$ &$\o{25}$  &$\o{145}$ &$\o{1245}$\\
$\o{13}\o{124}$ &$\o{13}\o{234}$ &$\o{124}\o{234}$ &$\o{125}$ &$\o{145}$  &$\o{245}$ \\
$\o{23}\o{14}$ &$\o{23}\o{1234}$ &$\o{14}\o{1234}$ &$\o{15}$  &$\o{245}$ &$\o{1245}$\\
$\o{23}\o{24}$ &$\o{23}\o{34}$ &$\o{24}\o{34}$     &$\o{25}$  &$\o{145}$ &$\o{1345}$ \\
$\o{23}\o{124}$ &$\o{23}\o{134}$ &$\o{124}\o{134}$ &$\o{125}$ &$\o{145}$  &$\o{245}$\\
$\o{123}\o{14}$ &$\o{123}\o{234}$ &$\o{14}\o{234}$ &$\o{15}$  &$\o{245}$ &$\o{1245}$ \\
$\o{123}\o{24}$ &$\o{123}\o{134}$ &$\o{24}\o{134}$ &$\o{25}$  &$\o{145}$ &$\o{1245}$\\
$\o{123}\o{124}$ &$\o{123}\o{34}$ &$\o{124}\o{34}$ &$\o{125}$ &$\o{145}$  &$\o{245}$\\ \hline
\end{tabular}\tae
Five of the rows can be derived by rotating the $bc$ indices positively and one negatively, wrapping $1\rightarrow4$. The other two rows, $\o2\o{13} \rightarrow \o{13}\o4$ and $\o{12}\o{13} \rightarrow \o{123}\o4$, do not follow a pattern nor do the $d$ columns when compared with \tar{zd12o}, so such pattern recognition is ignored. Since \tar{zd12} contains AAA non-associative cycles this gives another $21\times4\times3 = 252$ zero divisor pairs.

Now all other $bc$ pairs that don't include $\o5$ can be covered. This involves triads with overlapping AAA and BBA cycles but the AAA cases still have four modes and the BBA have prime and extended modes only. \tar{sai} shows that for $\PB{12}$ and $\PB{14}$ the type B triads provide half of the 12 zero divisor pairs for the algebra. In \tar{zd8} the $d$ column provides the AAA structure and the $d_1$ and $d_2$ columns the BBA structure. The 21 sets of 3-cycles in \tar{zd8} thus represent $21\times3\times4=252$ zero divisor pairs in the $d$ column and the $d_1$ and $d_2$ columns provide $21\times3\times2=126$ solutions each which means the table provides $252+126\times2=504$ zero divisor pairs.

\tar{zd8} can be replicated, apart from the $d$ column, with $\o4$ replaced with $\o5$ in the $bc$ cycles, as shown in \tar{zd4}. The $d_1$ and $d_2$ columns have the same BBA structure and modes as \tar{zd8} which represents $21\times3\times4=252$ solutions. \tas{zd8} and \ref{tab:zd4} provide the only B non-associative triads in $d_1$ and $d_2$ columns giving $21\times3\times2\times2=252$ triads, as expected.

Adding all solutions from  \tas{zd12o} through \ref{tab:zd4} gives $252 +252 +504 +252 = 1260$ which is enumerated in \tar{zdall} under the $\UB2$ zero divisors. Counting the number elements in the $d$ columns of these tables identifies 147 prime zero divisor pairs that generate the rest via the Cycle and Mode Theorems taking into account that some modes generate duplicates for some silos.

\tab{ht}{zd4}{BBA only Primary Zero Divisors in $\UB2$}
\begin{tabular}{|c|c|c||c|c|c|}  \hline 
 {\bf\em bc} &{\bf\em d$_1$} &{\bf\em d$_2$} &
 {\bf\em bc} &{\bf\em d$_1$} &{\bf\em d$_2$}\\\hline
$\o{1}\o{25}$    &$\o{345}$ &$\o{2345}$ &
$\o{3}\o{125}$   &$\o{145}$ &$\o{245}$\\
$\o{1}\o{35}$    &$\o{245}$ &$\o{2345}$ &
$\o{13}\o{15}$   &$\o{245}$ &$\o{1245}$\\
$\o{1}\o{235}$   &$\o{245}$ &$\o{345}$ &
$\o{13}\o{25}$   &$\o{145}$ &$\o{245}$\\
$\o{2}\o{15}$    &$\o{345}$ &$\o{1345}$ &
$\o{13}\o{125}$  &$\o{145}$ &$\o{245}$\\
$\o{2}\o{35}$    &$\o{145}$ &$\o{1345}$ &
$\o{23}\o{15}$   &$\o{245}$ &$\o{1245}$\\
$\o{2}\o{135}$   &$\o{145}$ &$\o{345}$ &
$\o{23}\o{25}$   &$\o{145}$ &$\o{1245}$\\
$\o{12}\o{15}$   &$\o{345}$ &$\o{1345}$ &
$\o{23}\o{125}$  &$\o{145}$ &$\o{245}$\\
$\o{12}\o{35}$   &$\o{145}$ &$\o{1345}$ &
$\o{123}\o{15}$  &$\o{245}$ &$\o{1245}$\\
$\o{12}\o{135}$  &$\o{145}$ &$\o{345}$ &
$\o{123}\o{25}$  &$\o{145}$ &$\o{1245}$\\
$\o{3}\o{15}$    &$\o{245}$ &$\o{2345}$ &
$\o{123}\o{125}$ &$\o{145}$ &$\o{245}$\\
$\o{3}\o{25}$    &$\o{145}$ &$\o{1245}$ &&&\\
 \hline
\end{tabular}\tae

\tar{zdall} provides the number of non-associative elements and zero divisors for the ultracomplex algebras up to level 10. These are labelled $\UB{m}$ where $n=m+3$ is the number of generators. Each generator level defines an algebra with $N=2^n-1$ pure basis elements and $\binom{N}3$ ordered product triples. Removing the associative triads and dividing by 28 provides the number of octonion and quasi-octonion algebra copies which corresponds with the multiples of 84 zero divisor pairs in the next level ultracomplex algebra. The zero divisors for $\UB1$ correspond to the AAA silo and those in $\UB2$ to all of the triads in the AAA and BBA silos. At higher levels, the silos after the first two have overlapping zero divisor pairs. Algebra $\UB3$ has new zero divisors from the ACC and XBB silos that need to be filtered but does not require the BXC silo. Algebras $\UB4$, $\UB5$, $\UB6$ and $\UB7$ add new zero divisors from the BXC silo and non-cycle triads but not the CAB silo.
\tab{ht}{zdall}{Ultracomplex Algebra Cardinality}
\begin{tabular}{|c|c|c|c|c|c|}  \hline   
\multirow{2}{*}{\bf Label} &{\bf Pure Basis} &\multicolumn{2}{|c|}{\bf Non-Associative}  &\multicolumn{2}{|c|}{\bf Zero Divisors}  \\\cline{3-6}
 &{\bf Size} &{\bf Triads} &{\bf 28 Factor} &{\bf Count} &{\bf 84 Factor} \\\hline
$\OB$ &7 &28  &1 &0 &0 \\ 
${\widetilde\OB}$ &7 &28 &-- &12 &-- \\
$\UB1\equiv\SB$ &15 &420 &15 &84 &1 \\ 
$\AA{3,1}$ &15 &420 &-- &84  &-- \\
$\AA{0,4}$ &15 &420 &-- &180 &-- \\
$\UB2$     &31   &4,340 &155 &1,260 &15 \\ 
$\UB3$     &63   &39,060 &1,395  &13,020 &155 \\ 
$\UB4$     &127  &330,708 &11,811 &117,180 &1,395 \\ 
$\UB5$     &255  &2,720,340 &97,155 &992,124 &11,811 \\ 
$\UB6$     &511  &22,064,980 &788,035 &8,161,020 &97,155 \\ 
$\UB7$  &1023  &177,736,020 &6,347,715 &66,194,940 &788,035 \\ 
\hline
\end{tabular}\tae

\section{Split Cayley-Dickson Algebras}
The split algebras are denoted $\AA{q,p}$ with $q$ being the number of negative and $p>1$ being positive signature generators. Hence $\AA{n}=\AA{n,0}$. Positive signature elements, called unitary generators, are denoted $\u{i}$ with $\epsilon_i=-1$ in (\ref{eqn:cd1}), $i\in\N_1^p$, so that $\u{i}^2=+1$, and split complex numbers (also known as Lorentz numbers), split quaternions and split octonions are denoted $\SCB$, $\SHB$ and $\SOB$, respectively~\cite{Harvey}. Similarly to the imaginary case, unitary multiplication, $\u\alpha\u\beta=\u\delta$ is also defined via the XOR operation $\delta=\alpha\veebar\beta$. The sign of this product can be based on the equivalent imaginary product except that cancellation of a unitary index does not negate the result. Mixtures of imaginary and unitary basis elements anticommute, $\o\alpha\u\beta=-\u\beta\o\alpha$, and non-associativity applies by assuming all unitary indices are greater than the imaginary indices. Lemma~\ref{lem:3} shows that squaring $\u\beta$ is $-1$ for even length of $\beta$ and is $+1$ for an odd length. Due to the symmetry of Cayley-Dickson algebras the assignment of indices to $\alpha$ and $\beta$ is irrelevant, it is only the number of indices that is relevant. This is now applied to each level of the graded algebras. 

The binomial expansion theorem is enumerated by Pascal's triangle that starts as\par
{\hfill
\begin{tabular}{ccccccccc}
  &  &  &  &1 &  &  &  &  \\
  &  &  &1 &  &1 &  &  &  \\
  &  &1 &  &2 &  &1 &  &  \\
  &1 &  &3 &  &3 &  &1 &  \\
1 &  &4 &  &6 &  &4 &  &1 \\
\end{tabular}\hfill}\vspace{2ex}\break
and is described by the recurrence relation that each position after the first is the sum of the two elements above. For graded algebras the triangle enumerates Cayley-Dickson algebras starting with the reals, $\R$.
Each row after this is the statement that the number of terms of any $k$-grade is given by the number of terms at grade $k$ without the generator $\o{n}$ or $\u{n}$ plus those that include this generator from grade $k-1$. Pascal's triangle is horizontally symmetric because for each term of grade $k$, the complementary term with grade $n-k$ contains the same $k$ indices removed from the $n$ generator case. For $n$ odd there are an equal number of terms with $\o{n}$ or $\u{n}$ as without in the complementary terms. For $n$ even the same applies apart from the middle grade $k=n/2$. But $\binom{n}{n/2}$ is always even for $n$ even and by the definition of combinations it distributes the $n$ generators equally into an even set, half of which contain $\o{n}$ or $\u{n}$ and the others exclude it. Of course, the left-hand side has grade $k=0$ that specifies the basis number 1 while the right-hand side, the complement of 1, has $k=n$ and always contains $\o{n}$ or $\u{n}$.

\begin{definition} The {\it pure trace} of the multiplication table for the algebra with $n$ generators is the sum of the squares of the $2^n-1$ pure basis terms, which is $-(2^n-1)$ for Cayley-Dickson algebra at level $n$. Thus for octonions this is $-7$.
\end{definition}

\begin{theorem}[Trace theorem]\label{thm:tt}
The pure trace for all of the split ultracomplex algebras is one, as it is for $\SCB$, $\SHB$ and $\SOB$~\cite{Harvey}.
\end{theorem}
\begin{proof} For any blade $\o\alpha\u\beta$, $\alpha$ and $\beta$ from disjoint sets of $\MB{n}$ then Lemma~\ref{lem:3} shows that $(\o\alpha\u\beta)^2=\pm1$, $-1$ if $\beta$ has even length and $+1$ if $\beta$ has odd length. Starting at $n=1$ the Lorentz numbers, $\SCB$, have one blade, $\u1$, so has pure trace one. Split quaternions, $\SHB$, can have pure basis, $\{\o1, \u2, \o1\u2\}$ or $\{\u1, \u2, \u{12}\}$, both of which satisfy the theorem. By induction, the basis at level $n-1$ is assumed to have pure trace one and either $\o{n}$ or $\u{n}$ will multiplying all basis elements. By the binomial theorem recurrence relation the existing pure basis is doubled by multiplying all existing basis elements by $\o{n}$ and adding the $\o{n}$ (multiplication by 1). Extending each $\alpha$ in this way does not change the parity of $\beta$ so this doubles the pure trace count and adds the -1 trace part for $\o{n}$, keeping the pure trace to one.

Now applying the same recurrence relation with $\u{n}$ swaps even and odd parities of each $\beta$ in the new multiplied terms when excluding the unit, 1. Hence a trace of one at level $n-1$ has terms added with opposite trace parity giving a total parity of zero. But the extra term, $\u{n}$, from multiplication by 1, regains the pure trace of one.
\end{proof}

The unitary Class~1, Class~2 and Class~3 associativity scheme for triads does not change from the non-split case because the same sign changes apply to each associativity term. This means the non-associativity types A, B, C and X are also appropriate but the zero divisors do change. These start appearing for the split octonions, $\SOB$, which have $12$ zero divisor pairs that, of course, have X-non-associativity triads. These are second and third cycles for a first cycle that is not a zero divisor and have extended dual zero divisor modes so only six need to be shown,
\begin{align*}
  &(\u1, \u{12}, \u3), (\u1, \u{13}, \u{123}), (\u2, \u{23}, \u{123}), \\
  &(\u2, \u{12}, \u3), (\u3, \u{13}, \u{123}), (\u3, \u{23}, \u{123})\es
\end{align*}
Split octonions have idempotents of the form $(1+\u{\alpha})/2$ with zero divisors $(1+\u{\alpha})(1-\u{\alpha})$ for $\alpha\in\N_1^3$ or $\alpha=\{123\}$. These are derived from the four positive square blades of $\SOB$, $\u{\alpha}^2=1$ plus $\u{123}^2=1$. Each idempotent can be multiplied on the left by the three blades that do not contain $\u{\alpha}$, which form a left ideal and gives $16$ zero divisors. Such cases are not covered by Theorem~\ref{thm:zd} that assumed non-scalar elements and multiplying by different terms can generate a nilpotent. For example, $(1+\u1)$ has $(-\u2+\u{12})(\u3+\u{13})=0$, which satisfies the definition for a zero divisor since $\u{12}\u3\u{13}=-\u2$ but  $(-\u2+\u{12})^2=0$. Hence nilpotents overlap with zero divisor pairs for split algebras.

Each Cayley-Dickson algebra, $\AA{q,p}$, at level $n=q+p$, is a doubling of a previous algebra at level $n-1$ by an extension of the basis with a new generator, $\gamma=\o{q+1}$ or $\gamma=\u{p+1}$. \tar{spa} shows the basis of the four split octonions and the extension to split sedenions by increasing $p$.
\tab{ht}{spa}{Split Sedenions}
\begin{tabular}{lll}
{\bf $\AA{q,p}$, $p+q=3$} & {\bf \qquad\qquad Basis} & {\bf $\AA{q,p}$, $p+q=4$, $p>0$} \\ \hline
$\AA{3}\cong\OB$  &$\{1, \o1, \o2, \o{12}, \o3, \o{13}, \o{23}, \o{123}\}$ &$\AA{3,1}\cong\AA{3}\oplus\AA{3}\u1$\\
$\AA{2,1}$  &$\{1, \o1, \o2, \o{12}, \u1, \o1\u1, \o2\u1, \o{12}\u1\}$   &$\AA{2,2}\cong\AA{2,1}\oplus\AA{2,1}\u2$\\
$\AA{1,2}$  &$\{1, \o1, \u1, \o1\u1, \u2, \o1\u2, \u{12}, \o1\u{12}\}$ &$\AA{1,3}\cong\AA{1,2}\oplus\AA{1,2}\u3$\\
$\AA{0,3}$  &$\{1, \u1, \u2, \u{12}, \u3, \u{13}, \u{23}, \u{123}\}$     &$\AA{0,4}\cong\AA{0,3}\oplus\AA{0,3}\u4$\\
\end{tabular}\tae

It is easy to see that the first two rows have structure $\HB+\HB\gamma$, for $n=3$, where $\gamma=\o3$ and $\u1$, respectively, and each split octonion in the last two rows have ${\widetilde\HB}+{\widetilde\HB}\gamma$, where $\gamma = \u2$ and $\u3$, respectively.
The middle two rows can also be interpreted as ${\widetilde\HB}+{\widetilde\HB}\gamma$, where $\gamma = \o2$ and $\o1$, respectively.
The proof that there are only two algebras satisfying $ab=\pm c$, for the three squares being $\pm 1$ is simply that $(ab)^2 = -baab = -a^2b^2=c^2$. Hence $a^2=b^2$ or one of these equals $c^2$, which means the signature is either all $-1$ or two $+1$ and one $-1$. Reorganising terms so that $c^2=-1$ and swapping $a$ and $b$ if necessary so that $abc=-1$ then the two split-quaternions $\{\u1,\u2,\u{12}\}$ and $\{\o1\u1,\u1,\o1\}$ are found to have identical structure and are thus isomorphic so that multiplying both of these by $\o2$ keeps the structure the same. This validates the well known result that the split-octonions are all isomorphic to each other~\cite{Harvey} and since $\AA{2,2}$, $\AA{1,3}$ and $\AA{0,4}$ are generated with a equivalent unitary generator then these are also isomorphic to each other. This is an example of Brown's more general isomorphism theorem~\cite{Brown} since $\AA{2,1}$, $\AA{1,2}$ and $\AA{0,3}$ are isomorphic and $\u2$, $\u3$, and $\u4$ have the same signature. For the remainder $\AA{0,3}$ is used for the ${\widetilde\OB}$ basis and $\AA{0,4}$ is used for the last three rows of split sedenions in \tar{spa}.

\tab{ht}{sszd}{Primary Zero Divisors for $\AA{0,4}$}
\begin{tabular}{|c|c|c|c||c|c|c|c|}  \hline   
{\bf\em b} &{\bf\em c} &{\bf\em d} &{\bf Modes} &{\bf\em b} &{\bf\em c} &{\bf\em d} & {\bf Modes} \\\hline
  $\u{1}$ &$\u{2}$   &$\u{13}, \u{14}, -\u{134}$ &2,2,2 &
  $\u{12}$ &$\u{3}$  &$\u{34}, -\u{134}$         &2,2 \\
  $\u{1}$ &$\u{3}$   &$\u{23}, \u{14}, -\u{124}$ &1,2,2 &
  $\u{12}$ &$\u{13}$ &$-\u{4}, \u{14}$           &4,4 \\
  $\u{1}$ &$\u{23}$  &$-\u{4}, \u{24}$           &2,2   &
  $\u{12}$ &$\u{4}$  &$\u{34}$                   &2 \\
  $\u{1}$ &$\u{4}$   &$\u{24}, \u{34}, \u{1234}$ &1,1,1 &
  $\u{12}$ &$\u{14}$ &$-\u{134}$                 &2 \\
  $\u{1}$ &$\u{24}$  &$-\u{134}, -\u{234}$       &1,1  &
  $\u{3}$  &$\u{4}$  &$\u{14}, \u{24}$           &1,1 \\
  $\u{1}$ &$\u{34}$  &$-\u{234}$                 &1   &
  $\u{3}$  &$\u{14}$ &$-\u{124}$                 &1 \\
  $\u{2}$ &$\u{3}$   &$\u{13}, \u{24}, -\u{124}$ &1,2,2 &
  $\u{3}$  &$\u{24}$ &$-\u{124}$                 &1 \\
  $\u{2}$ &$\u{13}$  &$\u{4}, \u{14}$            &2,2  &
  $\u{13}$ &$\u{4}$  &$\u{24}$                   &1 \\
  $\u{2}$ &$\u{4}$   &$\u{14}, \u{34}, \u{1234}$ &1,1,1 &
  $\u{13}$ &$\u{14}$ &$-\u{124}$                 &1 \\
  $\u{2}$ &$\u{14}$  &$-\u{134}, -\u{234}$       &1,1 &
  $\u{23}$ &$\u{4}$  &$\u{14}$                   &1 \\
  $\u{2}$ &$\u{34}$  &$-\u{134}$                 &1   &
  $\u{23}$ &$\u{24}$ &$-\u{124}$                 &1 \\ \hline
\end{tabular}\tae
Split sedenions, $\AA{0,4}$, have 180 zero divisor pairs, shown in \tar{sszd}, derived from 96 X, 48 A and 36 C-non-associative triads which includes extending zero divisors to allow for unitary products. Unlike the split octonions, $\AA{0,4}$ does have full cyclic structure but the second and third cycles sometimes require a different mode. \tar{sszd} documents the number of modes applicable for the $b,c,d$ triads in $d$ column order. Mode number 1 means only `p' is available, 2 means `p' and `d' or `e' is applicable, while 4 means all modes including `b'. The table shows that there are 39 primaries and 21 excess modes for these $d$ elements, which gives the $(39+21)\times3=180$ zero divisor pairs for $\AA{0,4}$. The $d$ blades shown as negative in the table require a variation of the zero divisor formula for some split triads, $(-a+b)(c+d)$. Taking the negative signature part of $\AA{0,4}$, $\AA{0,4}^-=\{\u{12}, \u{13}, \u{23}, \u{14}, \u{24}, \u{34}, \u{1234}\}$ then the multiplication table built from this basis is not $\OB$ but $\PB4$, which explains twelve of the zero divisor pairs.

The other split sedenion algebra in \tar{spa}, $\AA{3,1}=\OB\oplus\OB\u1$, has structure similar to $\AA4=\OB\oplus\OB\o4$. These have multiplication tables of the following form with vector $\vv{b}$ being the pure basis of $\OB$ and $\vv{b}^T$ being it's transpose and $\OB^2$ is the multiplication table for octonions and $\bar\OB^2$ is $\OB^2$ with positive signature or $+1$'s on the diagonal,
\begin{center}
\begin{tabular}{c|c|c|c|}
  \multicolumn{4}{c}{\qquad\quad\ $\vv{b}^T$\qquad$\gamma$\qquad$(\vv{b}\gamma)^T$} \\\cline{2-4}
 $\vv{b}\tstrut$    &$\OB^2$        &$\vv{b}\gamma$ &$-\bar\OB^2 \gamma$ \\[1ex]\cline{2-4}
 $\gamma$\tstrut    &$-\vv{b}^T\gamma$ &$\pm1$      &$\vv{b}$     \\[1ex]\cline{2-4}
 $\vv{b}\gamma$\tstrut &$-\OB^2\gamma$ &$\vv{b}$    &$\pm\OB^2$\\[1ex]\cline{2-4}
\end{tabular}
\end{center}
where $\gamma = \u4$ and $+$ on the diagonal for $\AA{3,1}$ and $\gamma=\o4$ and $-$ for $\AA{4,0}$. These algebras both have 84 zero divisor pairs, as shown for $\AA{4}$ in \tar{zds} and the $\AA{3,1}$ zero divisors have analogous structure with seven primary pairs that expand to three cyclic pairs and four modes as shown in \tar{zd31}. In this case all modes have $d=\u1$. Thus it is conjectured that this split algebra involves seven quasi-octonions with mixed signature that are a variation of $\PB4$. This work verifies, for sedenions, the work of Brown~\cite{Brown} and Eakin and Sathaye~\cite{Eakin} showing that the Cayley-Dickson algebra structure depends on the source algebra and the basis multiplier, $\gamma$.

\tab{ht}{zd31}{Primary Zero Divisors for $\AA{3,1}$}
\begin{tabular}{|c|c|c|c|c|c|}  \hline 
  $b$    &$c$      &$d$    &$a$           &$db$         &$ab$\\\hline
  $\o1$  &$\o2$    &$\u1$  &$\o{12}\u1$   &$-\o1\u1$    &$-\o2\u1$ \\
  $\o1$  &$\o3$    &$\u1$  &$\o{13}\u1$   &$-\o1\u1$    &$-\o3\u1$ \\
  $\o1$  &$\o{23}$ &$\u1$  &$-\o{123}\u1$ &$-\o1\u1$    &$-\o{23}\u1$ \\
  $\o2$  &$\o3$    &$\u1$  &$\o{23}\u1$   &$-\o2\u1$    &$-\o3\u1$ \\
  $\o2$  &$\o{13}$ &$\u1$  &$\o{123}\u1$  &$-\o2\u1$    &$-\o{13}\u1$ \\
  $\o{12}$ &$\o3$ &$\u1$  &$\o{123}\u1$  &$-\o{12}\u1$ &$-\o3\u1$ \\
  $\o{12}$ &$\o{13}$ &$\u1$ &$-\o{23}\u1$ &$-\o{12}\u1$ &$-\o{13}\u1$ \\\hline
\end{tabular}\tae

An overlap between split and non-split zero divisor pairs does not occur. The triad $(\u{12},\u{13},\u{14})$ and its product, $\u{12}\u{13}\u{14}=-\u{1234}$, contains terms that are all imaginary, which should be analogous to the non-split triad $(\o{12},\o{13},\o{14})$ but this has product $+\o{1234}$ providing a zero divisor pair with opposite parity.

\section{Summary}
The Cayley-Dickson algebras have been assumed to be symmetric in the permutation of basis elements but this is only valid for hypercomplex numbers with octonions introducing non-associativity but in a symmetric, alternating form. Using the graded notation provides a rigorous definition to the ordering of basis elements and exposes the asymmetry of the power-associative Cayley-Dickson algebras. The introduction of the terminology {\it ultracomplex numbers} allows power associative Cayley-Dickson and the quasi-Cayley-Dickson subalgebras to be analysed and found to be more complex than hypercomplex numbers, justifying the new terminology. It is not just the Cayley-Dickson algebras that can be power-associative and non-alternative but their subalgebras, such as the three quasi-octonions and the twelve quasi-sedenions that can be found inside $\UB2$. The most important step was the introduction of the ordered associator to identify the different types of non-associativity thus allowing the full range of power-associative subalgebras of Cayley-Dickson algebras to be classified. This was enabled by the grading based on binary expansion to uniquely identify the multiples of 84 zero divisor pairs of two terms. Developing cycles and modes to analyse the algebra structure, the non-associative rings can be reduced to just seven triads that are the prime elements of the quasi-octonions. The total number of zero divisor pairs for each Cayley-Dickson algebra is also derived. 

The graded approach also applies to split ultracomplex numbers by finding the isomorphisms between some of the split sedenions. It is convenient to use $\epsilon=-1$ for all generators for the isomorphic cases, which then have positive square generators, $\u{i}^2=1$. The convenient symmetry thereafter is that all odd length grades are unitary and all even grades are imaginary. It would be interesting to extend these results to the three split $\UB2$ algebras predicted by Brown and Eakin and Sathaye.

The non-associativity from the non-split cases carry over to split algebras for equivalent indices but some of the zero divisors are related to but not the same as ideals of the algebra. Open questions for future investigation is to establish connections of split zero divisors to ideals, gain further understanding of how zero divisors relate to Malcev algebras and investigate the relationships between the quasi-octonion algebras to establish why only three of the six are contained in power-associative algebras. The use of the calibration from Clifford algebra~\cite{Wilmot} for generating all Cayley-Dickson algebras also needs to be formalised.

\ifnum1=\ELS
The author gratefully acknowledges the invaluable assistance and encouragement of James Chappell, Jacqui Ramagge and Derek Abbott for improving terminology, mathematical understanding and ongoing assistance with the publishing of this article. Acknowledgement is also due to Maher ben Abdessalem, Esprit, Tunisia for communicating a correction to $Z_m$ in an earlier draft of the manuscript.

All of this work was verified with the use of an octonion/Cayley-Dickson calculator written in Python. The code is {\it calcO.py} using the {\it load(ultronions)} module. Specifically, the {\it triadCycles} function generates the data for the Cycles Theorem. The github URL for the calculator is\par\centerline{\url{https://github.com/GPWilmot/geoalg}.}
\bibliography{Article2} 
\else
\ifnum1=\RSOS
All of this work was verified with the use of an octonion/Cayley-Dickson calculator written in Python. The code is {\it calcO.py} using the {\it load(ultronions)} module. Specifically, the {\it triadCycles} function generates the data for the Cycles Theorem. The github URL for the calculator is\par\centerline{\url{https://github.com/GPWilmot/geoalg}.}
\funding{The author gratefully acknowledges the support of an Australian Government Research Training Program Scholarship and from the Australian Research Council (FL240100217).}
\ack{The author acknowledges the invaluable assistance and encouragement of James Chappell, Jacqui Ramagge and Derek Abbott for improving terminology, mathematical understanding and ongoing assistance with the publishing of this article. Acknowledgement is also due to Maher ben abdessalem, Espirit, Tunisia for communicating a correction to $Z_m$.}

\else
The author acknowledges the invaluable assistance and encouragement of James Chappell, Jacqui Ramagge and Derek Abbott for improving terminology, mathematical understanding and ongoing assistance with the publishing of this article. Acknowledgement is also due to Maher ben abdessalem, Espirit, Tunisia for communicating a correction to $Z_m$.
The author also gratefully acknowledges the support of an Australian Government Research Training Program Scholarship and from the Australian Research Council (FL240100217).
All of this work was verified with the use of an octonion/Cayley-Dickson calculator written in Python. The code is {\it calcO.py} using the {\it load(ultronions)} module. Specifically, the {\it triadCycles} function generates the data for the Cycles Theorem. The github URL for the calculator is\par\centerline{\url{https://github.com/GPWilmot/geoalg}.}

\fi\fi
\end{document}